\documentclass[12pt]{article}
\usepackage{graphicx}
\usepackage{natbib}
\usepackage{enumerate}
\usepackage{color}

\usepackage{amssymb}
\usepackage{amsmath}
\usepackage{amsfonts}
\usepackage{amsthm}

\newcommand{\tr}{^{\mbox{\scriptsize tr}}}

\newcommand{\R}{\mbox{\bf R}}
\newcommand{\C}{\mbox{\bf C}}

\newcommand{\E}{\mbox{\bf E}}

\newcommand{\var}{\mbox{\bf Var}}

\newcommand{\e}{\mbox{\rm e}}

\newcommand{\bec}{\begin{equation}}
\newcommand{\eec}{\end{equation}}

\newcommand{\bm}[1]{{\mbox{\boldmath $#1$}}}
\newcommand{\sbm}[1]{{\mbox{\scriptsize\boldmath $#1$}}}
\newcommand{\samp}[2]{\ifmmode #1_1,\dots,#1_{#2}\else
$#1_1,\dots,#1_{#2}$\fi}
\newcommand{\tras}{^{\mbox{\scriptsize tr}}}

\newtheorem{theorem}{Theorem}

\newtheorem{lemma}{Lemma}

\newtheorem{corollary}{Corollary}

\newcommand{\pf}{{\noindent {\em Proof: }}}

\theoremstyle{definition}

\newtheorem{example}{Example}

\begin{document}

\bibliographystyle{plainnat}

\title{Modeling stationary data by a class of generalised Ornstein-Uhlenbeck processes.}
\author{Argimiro Arratia\thanks{Dept.  Llenguatges i Sistemes Inform\`atics, Universitat Polit\`ecnica de Catalunya, Barcelona, Spain. 
Supported  by BASMATI MICINN project (TIN2011-27479-C04-03), 
SGR2009-1428 (LARCA) 
 and SINGACOM (MTM2007-64007)}\hspace{1mm}, Alejandra Caba\~na\thanks{Departament de Matem\`atiques, Universitat Aut\`onoma de Barcelona, Spain. Partially supported by TIN2008-06582-C03-02, Ministerio de Ciencia y Tecnolog\'{\i}a, Spain}\hspace{1mm}
  and
Enrique M. Caba\~na \thanks{Departamento de M\'etodos Matem\'atico-Cuantitativos,  Universidad de la
Rep\'ublica, Montevideo, Uruguay.}}

\date{}

\maketitle

\begin{abstract}
An Ornstein-Uhlenbeck (OU) process can be considered as a continuous time interpolation of 
the discrete time AR$(1)$ process. Departing from this fact, we analyse  in this work the effect
of iterating OU treated as a linear operator that maps a Wiener process onto Ornstein-Uhlenbeck process, so as to build a family of higher order Ornstein-Uhlenbeck processes, OU$(p)$, in a similar spirit as the higher order autoregressive processes AR$(p)$. 
We show that for $p \ge 2$ we obtain in general a process with  covariances different than  those of an AR$(p)$, and that for various continuous time processes, sampled from real data at equally spaced time instants, the OU$(p)$ model outperforms the appropriate AR$(p)$ model.
 Technically our composition of the OU operator is easy to manipulate and its parameters can be computed efficiently because, as we show, the iteration of OU operators leads to a process that can be expressed as a linear combination of basic OU processes. Using this expression we obtain a closed formula for the covariance of the iterated OU process, and consequently estimate the parameters of an OU$(p)$ process by maximum likelihood or, as an  alternative, by matching correlations, the latter being  a procedure resembling the method of moments.

\bigskip


\end{abstract}

\section{Introduction}

The Ornstein-Uhlenbeck process  (from now on OU) 
was introduced by L. S. Ornstein and E. G. Uhlenbeck \citep{OU} as a model for the velocities of a particle subject to the collisions   with surrounding molecules. It improves Einstein's model (a Wiener process) because it also applies to fluids with finite viscosity, and since the 1950's is a well 
studied  and accepted  model for thermodynamics, chemical and other various stochastic processes found in physics and the natural sciences \citep{Gard04}. 
Moreover, the OU process is the unique non--trivial stochastic process
that is stationary, Markovian and Gaussian \citep{Maller09}.
Additionally it is mean-reverting, and for all these properties it has found its way into 
financial engineering, first  as a model for the term structure of interest 
rates in a form due to \cite{Vas77}, and then under other variants or generalisations (e.g. where the underlying random noise is
a L\'evy process) as a model of financial time series with applications to 
option pricing, portfolio optimisation and risk theory, among others
\citep[and references there in]{NV2003,BS2001a,BS2001b,Maller09}.

The OU process can be thought of as  continuous time interpolation of an autoregressive process of order one (i.e. an AR(1) process), 
a link that we shall make evident  in Section  \ref{prelim}.
Beginning with this relation to the autoregressive  model, one can seek 
to define and analyse  the result of iterating the application of the operator 
that maps a Wiener process on a OU process.
Thus,
in Section \ref{oupsec} a new family of processes is introduced: for each positive integer $p$ the
{\sl  Ornstein-Uhlenbeck processes of order $p$}, denoted OU($p$), are defined and proposed as models for either stationary continuous time processes or the series obtained by observing these continuous processes at equally spaced instants.  The OU(1) processes are the ordinary Ornstein - Uhlenbeck processes.
While the series obtained by sampling OU(1) processes at equally spaced times are autoregressive of the same order, this property does not extend in general for $p>1$ as shown in Section \ref{oupisnotarp}. Hence, OU processes of higher order appear as a new model, competitive in a discrete time setting with higher order autoregressive  processes (AR or ARMA).
The estimation of the parameters of OU($p$) processes is attempted in Section \ref{proce}, and examples showing the comparison of the proposed methods for that estimation and the application of OU($p$) models to real data are provided in Section \ref{examp}.
Section \ref{conclus} contains our concluding remarks.

\section{Preliminaries}\label{prelim}
Let us call $w$ a standard Wiener process, that is, a Gaussian, centred process with independent  increments with variance $\E(w(t)-w(s))^2=|t-s|$.
We impose further (as usual) that  $w(0)=0$, but shall not limit the domain of the parameter to $\R^+$ and assume that $w(t)$ is defined for $t$ in $\R$.
Then, an Ornstein-Uhlenbeck process with parameters $\lambda >0 ,\sigma>0 $ can be written as
\bec\label{OU0}
\xi_{\lambda,\sigma}(t)=\sigma \int_{-\infty}^t\e^{-\lambda(t-s)}dw(s)
\eec
or, in differential form,
\bec\label{eq:OUdif}
d\xi_{\lambda,\sigma}(t)=-\lambda \xi_{\lambda,\sigma}\;dt+\sigma dw(t)
\eec

We may think of $\xi_{\lambda,\sigma}$ as the result of accumulating a random noise, with reversion to the mean (that we assume to be 0) of exponential decay with rate
 $\lambda$. 
 The magnitude of the noise is given by $\sigma$. 
 
A widely used class of models for discrete time stationary series are the autoregressive process of order $p$, AR$(p)$. They are obtained from a series $W_t$ of standard Gaussian independent  random variables as 
$$(1-\phi_1B-\phi_2B^2-\dots-\phi_pB^p) x_t= \prod_{j=1}^p (1+\lambda_j B) x_t=\sigma W_t$$
where $B$ is the backshift operator that carries $x_t$ into $x_{t-1}$.  Moreover, the innovations $W_t$ can be thought of as $W_t= w(t)-w(t-1)$ for a standard Wiener process $w$. The process $x$ is stationary if $|\lambda_j| <1$ for all $j=1,2,\dots,p$.

When the Ornstein-Uhlenbeck process $x$ is sampled at equally spaced times $\{i\tau: i=0,1,2,\dots,n\}$, $\tau >0$, 
the series $X_i=x(i\tau)$ obeys an autoregressive model of order 1, AR(1),  since
$$X_{i+1}=\sigma \int_{-\infty}^{(i+1)\tau}\e^{-\lambda((i+1)\tau-s)}dw(s)$$
$$=\sigma\e^{-\lambda\tau}\int_{-\infty}^{i\tau}\e^{-\lambda(i\tau-s)}dw(s)+\sigma\int_{i\tau}^{(i+1)\tau}\e^{-\lambda((i+1)\tau-s)}dw(s)=\e^{-\lambda\tau}X_i+Z_i,$$
where $Z_i=\sigma\int_{i\tau}^{(i+1)\tau}\e^{-\lambda((i+1)\tau-s)}dw(s)$ is a Gaussian innovation (independent of $\{w(t): t\leq i\tau\}$ and 
$\{x(t): t\leq i\tau\}$) with variance $$\sigma^2\int_{i\tau}^{(i+1)\tau}\e^{-2\lambda((i+1)\tau-s)}ds
=\sigma^2\int_{-\tau}^0\e^{2\lambda s}ds=\frac{\sigma^2}{2\lambda}(1-\e^{-2\lambda\tau}).$$

Hence, we can consider the OU process as continuous time interpolation of an AR(1) process.
Notice that 
both models are stationary.
As we show in Section \ref{oupisnotarp}, the result of iterating the operator that carries Wiener process into Ornstein-Uhlenbeck process is  not an interpolation of an autoregressive process.

\section{Ornstein-Uhlenbeck processes of order $p$}
\label{oupsec}

Let ${\cal OU}_{\lambda}$ be defined as the operator that maps $\sigma w$ onto $\xi_{\lambda,\sigma}(t)$, and also maps a differentiable process $y(t),t\in\R$ onto \bec\label{elou}{\cal OU}_{\lambda}y(t)=\int_{-\infty}^t\e^{-\lambda(t-s)}dy(s),\eec
when the integral converges.
The definition is extended to include complex processes, by replacing $\lambda$ by $\kappa=\lambda+\imath\mu$, $\lambda>0$, $\mu\in\R$ in (\ref{elou}). The set of complex numbers with positive real part is denoted by $\C^+$.

For $p\ge1$, the process
\bec\label{oup}
x={\cal OU}_{\bm\kappa}(\sigma w):={\cal OU}_{\kappa_1}{\cal OU}_{\kappa_2}\cdots {\cal OU}_{\kappa_p}(\sigma w)=\prod_{j=1}^p{\cal OU}_{\kappa_j}(\sigma w)
\eec
will be called {\em Ornstein-Uhlenbeck process of order $p$ with parameters $\bm\kappa=(\kappa_1,\dots,\kappa_p)\in (\C^+)^p$} and $\sigma>0$.  The composition $\prod_{j=1}^p{\cal OU}_{\kappa_j}$ is unambiguosly defined because the application of ${\cal OU}_{\kappa_j}$ operators is commutative as shown in Theorem \ref{teo1} (\ref{p1}) below.

For technical reasons, it is convenient to introduce the {\em Ornstein-Uhlenbeck operator ${\cal OU}_{\kappa}^{(h)}$ of degree $h$ with parameter $\kappa$} that maps $y$ onto
\bec\label{ouph}
{\cal OU}_{\kappa}^{(h)}(t)y(t)=\int_{-\infty}^t\e^{-\kappa(t-s)}\frac{(-\kappa(t-s))^h}{h!}dy(s)
\eec
and $\sigma w$ onto
\bec\label{defxih}
\xi_{\kappa,\sigma}^{(h)}(t)=\sigma \int_{-\infty}^t\e^{-\kappa(t-s)}\frac{(-\kappa(t-s))^h}{h!}dw(s)
\eec

\subsection{Properties}\label{s2}
The following statements summarize some properties of products (compositions) of the operators defined by (\ref{oup}) and (\ref{ouph}), and correspondingly, of the stationary centred Gaussian processes $\xi_{\kappa,\sigma}^{(h)}$, $h=0,1,2,\dots$. These processes will be called {\em Ornstein-Uhlenbeck processes of degree $h$}. 
In particular, the Ornstein-Uhlenbeck processes of degree zero  $\xi_{\kappa,\sigma}^{(0)}=\xi_{\kappa,\sigma}$ are the ordinary Ornstein-Uhlenbeck processes (\ref{OU0}).

\begin{theorem}\label{teo1}
$ $

\begin{enumerate}[(i)]

\item\label{p1}
When $\kappa_1\not=\kappa_2$, the product ${\cal OU}_{\kappa_2}{\cal OU}_{\kappa_1}$ can be computed as $$\frac{\kappa_1}{\kappa_1-\kappa_2}{\cal OU}_{\kappa_1}+\frac{\kappa_2}{\kappa_2-\kappa_1}{\cal OU}_{\kappa_2}$$ and is therefore commutative.
\item\label{p2}
The composition $\prod_{j=1}^p{\cal OU}_{\kappa_j}$ constructed with values of $\kappa_1,\dots,\kappa_p$ pairwise different, is equal to the linear combination
\bec\label{lin}\prod_{j=1}^p{\cal OU}_{\kappa_j}=\sum_{j=1}^pK_j(\kappa_1,\dots,\kappa_p){\cal OU}_{\kappa_j}, \eec
with coefficients
\bec\label{laK}K_j(\kappa_1,\dots,\kappa_p)=\frac{1}{\prod_{\kappa_l\not=\kappa_j}(1-\kappa_l/\kappa_j)}.\eec
\item\label{p3}
For $i=1,2,\dots$, ${\cal OU}_{\kappa}{\cal OU}_{\kappa}^{(i)}={\cal OU}_{\kappa}^{(i)}-\kappa{\cal OU}_{\kappa}^{(i+1)}$.
\item\label{p4}
For any positive integer $p$ the $p$-th power of the Ornstein-Uhlenbeck operator has the expansion
\bec\label{potp}{\cal OU}_{\kappa}^p=\sum_{j=0}^{p-1}{p-1\choose j}{\cal OU}_{\kappa}^{(j)}.\eec
\item\label{p5}
Let $\kappa_1,\dots,\kappa_q$ be pairwise different complex numbers with positive real parts, and $p_1,\dots,p_q$ positive integers, and let us denote by $\bm\kappa$ a complex vector in $(\C^+)^p$ with components $\kappa_h$ repeated $p_h$ times,  $p_h\geq 1$, $h=1,\dots,q$, $\sum_{h=1}^qp_h=p$. Then, with $K_h(\bm\kappa)$ defined by (\ref{laK}),
$$\prod_{h=1}^q{\cal OU}_{\kappa_h}^{p_h}=\sum_{h=1}^q\frac{1}{\prod_{l\not=h}(1-\kappa_l/\kappa_h)^{p_l}}{\cal OU}_{\kappa_h}^{p_h}=\sum_{h=1}^qK_h(\bm\kappa){\cal OU}_{\kappa_h}^{p_h}.$$
\end{enumerate}
\end{theorem}

\begin{corollary}\label{defx} 
The process $$x= {\cal OU}_{\bm\kappa}(\sigma w)=\prod_{h=1}^q{\cal OU}_{\kappa_h}^{p_h}(\sigma w)$$ can be expressed as the linear combination 
\bec\label{lali}x=\sum_{h=1}^qK_h(\bm\kappa)(1+\xi_{\kappa_h,\sigma})^{(p_h-1)},\quad
(1+\xi_{\kappa_h,\sigma})^{(p_h-1)}=\sum_{j=0}^{p_h-1}\textstyle {p_h-1\choose j}\xi_{\kappa_h,\sigma}^{(j)}
\eec
of the $p$ processes $\{\xi_{\kappa_h,\sigma}^{(j)}: h=1,\dots,q, j=0\dots,p_h-1\}$ (see (\ref{defxih})).

\end{corollary}

\begin{corollary}\label{rea} For real $\lambda,\mu$,  with $\lambda>0$, the product ${\cal OU}_{\lambda+\imath\mu} {\cal OU}_{\lambda-\imath\mu}$ is real, that is, applied to a real process produces a real image.
\end{corollary}

\noindent
{\em Proof of the Theorem and its corollaries}:

Parts $(\ref{p1})$ and $(\ref{p3})$ are obtained by direct computation of the integrals, $(\ref{p2})$ follows from $(\ref{p1})$ by finite induction, as well as $(\ref{p4})$ from $(\ref{p3})$.

From the continuity of the integrals with respect to the parameter $\kappa$, the power ${\cal OU}_{\kappa}^p$ satisfies
\bec\label{elim}{\cal OU}_{\kappa}^p=\lim_{\delta\downarrow0}\prod_{j=1}^p{\cal OU}_{\kappa+j\delta}=\lim_{\delta\downarrow0}\sum_{j=1}^pK_j'(\delta,\kappa,p){\cal OU}_{\kappa+j\delta}\eec with $$K_j'(\delta,\kappa,p)=\frac{1}{\prod_{1\leq l\leq p,l\not=j}(1-\frac{\kappa+l\delta}{\kappa+j\delta})}.$$

On the other hand, by $(\ref{p1})$, \bec\label{suma}\prod_{h=1}^q{\cal OU}_{\kappa_h}^{p_h}=\lim_{\sbm\delta\downarrow0}\prod_{h=1}^q\prod_{j=1}^{p_h}{\cal OU}_{\kappa_h+j\delta_h}=\lim_{\sbm\delta\downarrow0}\sum_{h=1}^q\sum_{j=1}^{p_h}K_{h,j}''(\bm\delta,\bm\kappa){\cal OU}_{\kappa_h+j\delta_h}\eec
where $\bm\delta=(\delta_1,\dots,\delta_q)$,
$$K_{h,j}''(\bm\delta,\bm\kappa)=\frac{1}{\prod_{1\leq h'\leq q,1\leq j'\leq p_h,(h',j')\not=(h,j)}(1-\frac{\kappa_{h'}+j'\delta_{h'}}{\kappa_h+j\delta_h})}
=K_{h,j}'''(\bm\delta,\bm\kappa)K_j'(\delta_h,\kappa_h,p_h),$$ and $$K_{h,j}'''(\bm\delta,\bm\kappa)=\frac{1}{\prod_{1\leq h'\leq q,h'\not=h}\prod_{j'=1}^{p_{h'}}(1-(\kappa_{h'}+j'\delta_{h'})/(\kappa_h+j\delta_h))}
\rightarrow K_h(\bm\kappa)\mbox{ as }\sbm\delta\downarrow0$$

For the $h$-th term in the right-hand side of (\ref{suma}), we compute
$$\lim_{\sbm\delta\downarrow0}\sum_{j=1}^{p_h}K_{h,j}''(\bm\delta,\bm\kappa){\cal OU}_{\kappa_h+j\delta_h}
=\lim_{\sbm\delta\downarrow0}\sum_{j=1}^{p_h}K_{h,j}'''(\bm\delta,\bm\kappa)K_j'(\delta_h,\kappa_h,p_h){\cal OU}_{\kappa_h+j\delta_h}$$

$$=\lim_{\sbm\delta\downarrow0}\sum_{j=1}^{p_h}(K_{h,j}'''(\bm\delta,\bm\kappa)-K_h(\bm\kappa))K_j'(\delta_h,\kappa_h,p_h){\cal OU}_{\kappa_h+j\delta_h}$$ $$+K_h(\bm\kappa)\lim_{\sbm\delta\downarrow0}\sum_{j=1}^{p_h}K_j'(\delta_h,\kappa_h,p_h){\cal OU}_{\kappa_h+j\delta_h}=K_h(\bm\kappa){\cal OU}_{\kappa_h}^{p_h}$$
because of (\ref{elim}), since, in addition, each term in the first sum tends to zero. This ends the verification of $(\ref{p5})$.

Corollary \ref{defx} is an immediate consequence of  $(\ref{p4})$
 and $(\ref{p5})$, and Corollary \ref{rea} follows by applying $(\ref{p1})$,
  to compute $${\cal OU}_{\lambda+\imath\mu} {\cal OU}_{\lambda-\imath\mu}=\frac{\lambda+\imath\mu}{2\imath \mu}{\cal OU}_{\lambda+\imath\mu}-\frac{\lambda-\imath\mu}{2\imath \mu}{\cal OU}_{\lambda-\imath\mu}$$
$$=\int_{-\infty}^t\e^{-\lambda(t-s)}\left[\textstyle\frac{\lambda+\imath\mu}{2\imath \mu}(\cos(\mu(t-s))+\imath\sin(\mu(t-s)))\right.$$
$$\left.-\textstyle\frac{\lambda-\imath\mu}{2\imath \mu}(\cos(\mu(t-s))-\imath\sin(\mu(t-s)))\right]dw(s)$$
$$=\int_{-\infty}^t\e^{-\lambda(t-s)}(\cos(\mu(t-s))+\textstyle\frac{\lambda}{\mu}\sin(\mu(t-s)))dw(s).$$\qed

\subsection{Computing the covariances}
The representation $$x= {\cal OU}_{\bm\kappa}(\sigma w)
=\sum_{h=1}^qK_h\sum_{j=1}^{p_h}{{p_h-1}\choose {j-1}}{\cal OU}_{\kappa_h}^{(j-1)}(\sigma w)$$
of $x$ as a linear combination of the processes $\xi_{\kappa_h,\sigma}^{(i)}={\cal OU}_{\kappa_h}^{(i)}(\sigma w)$ allows a direct computation of the covariances 
 $\gamma(t)=\E x(t)\bar x(0)$ through a closed formula, in terms of 
the covariances
$\gamma_{\kappa_1,\kappa_2,\sigma}^{(i_1,i_2)}(t)=\E \xi_{\kappa_1,\sigma}^{(i_1)}(t)\bar\xi_{\kappa_2,\sigma}^{(i_2)}(0)$:

\bec\label{cf}\gamma(t)\!\!=\!\!\sum_{h'=1}^q\sum_{i'=0}^{p_{h'}-1}\sum_{h''=1}^q\sum_{i''=0}^{p_{h''}-1}\!\!K_{h'}(\bm\kappa)\bar K_{h''}(\bm\kappa){{p_{h'}-1}\choose i'}{{p_{h''}-1}\choose i''}\gamma_{\kappa_{h'},\kappa_{h''},\sigma}^{(i',i'')}(t)
\eec
with
$$\gamma_{\kappa_1,\kappa_2,\sigma}^{(i_1,i_2)}(t)=
\sigma^2(-\kappa_1)^{i_1}(-\bar\kappa_2)^{i_2}\int_{-\infty}^0\e^{-\kappa_1(t-s)}\frac{(t-s)^{i_1}}{i_1!}\e^{-\bar\kappa_2(-s)}\frac{(-s)^{i_2}}{i_2!}ds$$
$$=
\sigma^2(-\kappa_1)^{i_1}(-\bar\kappa_2)^{i_2}\e^{-\kappa_1t}\sum_{j=0}^{i_1}{i_1\choose j}\frac{t^j}{i_1!i_2!}\int_{-\infty}^0\e^{(\kappa_1+\bar\kappa_2) s}(-s)^{i_1+i_2-j}ds$$
\bec\label{laga}=\frac{\sigma^2(-\kappa_1)^{i_1}(-\bar\kappa_2)^{i_2}\e^{-\kappa_1t}}{i_2!}\sum_{j=0}^{i_1}\frac{t^j(i_1+i_2-j)!}{j!(i_1-j)!(\kappa_1+\bar\kappa_2)^{(i_1+i_2-j+1)}}\eec

A  real expression  for the covariance when the imaginary  parameters appear as conjugate pairs is much more involved than this one, that contains complex terms.

\section{OU($p$)   is not an AR($p$)}\label{oupisnotarp}
The series of observations of an OU($p$) at equally spaced times is not an AR($p$) in general, for $p>1$.
Consider the autocorrelations of the time series obtained by evaluating $x$ in multiples of a given instant
$\tau$. 
If $X_i=x(i\tau)$, then  $$\rho_i=\frac{\E X_iX_0}{\var X_0}=\frac{\gamma(i\tau)}{\gamma(0)}$$
Lemma \ref{leOU2AR}  shows through a direct computation 
of covariances that the family
of series obtained from OU(2) are not in general AR(2) processes.


\begin{lemma}\label{leOU2AR}
Suppose $x$ is a real OU(2) process and $X$ is a AR(2) process, 
 with equal
autocorrelations of order 1 and 2. Then, in general, the autocorrelations of order 3 are different.
\end{lemma}
\pf 
For the  AR(2) process, if $r_1, r_2, r_3$ are the autocorrelations 
of orders 1, 2 and 3, these quantities satisfy the following relations:
 \bec\label{lar3}-1\leq r_1\leq 1,\;\; 2r_1^2-1\leq r_2\leq 1, \;\;r_3=\frac{r_1}{1-r_1^2}(2r_2-r_1^2-r_2^2)
 \eec

For the  OU(2) process of real parameters $\lambda_1<\lambda_2$,
$$\gamma(t)=\sum_{j,k=1}^2\frac{K_jK_k\e^{-\lambda_j t}}{\lambda_j+\lambda_k}
=\frac{K_1^2\e^{-\lambda_1 t}}{2\lambda_1}+\frac{K_2^2\e^{-\lambda_2 t}}{2\lambda_2}+
\frac{K_1K_2(\e^{-\lambda_1 t}+\e^{-\lambda_2 t})}{\lambda_1+\lambda_1}$$
$$=\frac{\lambda_1\e^{-\lambda_1 t}}{2(\lambda_2-\lambda_1)^2}+\frac{\lambda_2\e^{-\lambda_2 t}}{2(\lambda_2-\lambda_1)^2}-\frac{\lambda_1\lambda_2(\e^{-\lambda_1 t}+\e^{-\lambda_2 t})}{(\lambda_2-\lambda_1)^2(\lambda_1+\lambda_2)}$$
$$=\frac{\lambda_1(\lambda_1+\lambda_2)\e^{-\lambda_1 t}+\lambda_2(\lambda_1+\lambda_2)\e^{-\lambda_2 t}-2\lambda_1\lambda_2(\e^{-\lambda_1 t}+\e^{-\lambda_2 t})}{2(\lambda_2-\lambda_1)^2(\lambda_1+\lambda_2)}$$
$$=\frac{\lambda_1(\lambda_1-\lambda_2)\e^{-\lambda_1 t}+\lambda_2(-\lambda_1+\lambda_2)\e^{-\lambda_2 t}}{2(\lambda_2-\lambda_1)^2(\lambda_1+\lambda_2)}
=\frac{\lambda_2\e^{-\lambda_2 t}-\lambda_1\e^{-\lambda_1 t}}{2(\lambda_2^2-\lambda_1^2)}$$ 
From this equation we obtain
$$\gamma(0)=\frac{1}{2(\lambda_1+\lambda_2)}$$
and the correlations 
\bec\label{larho}\rho_h=\frac{\lambda_2\e^{-\lambda_2h\tau}-\lambda_1\e^{-\lambda_1h\tau}}{\lambda_2-\lambda_1}, \;\; h=1,2,3,\dots\eec

 \begin{figure}[!t]

\vspace{-.8cm}

\noindent\hspace{-1cm}\raisebox{-2cm}{\includegraphics[width=9cm]{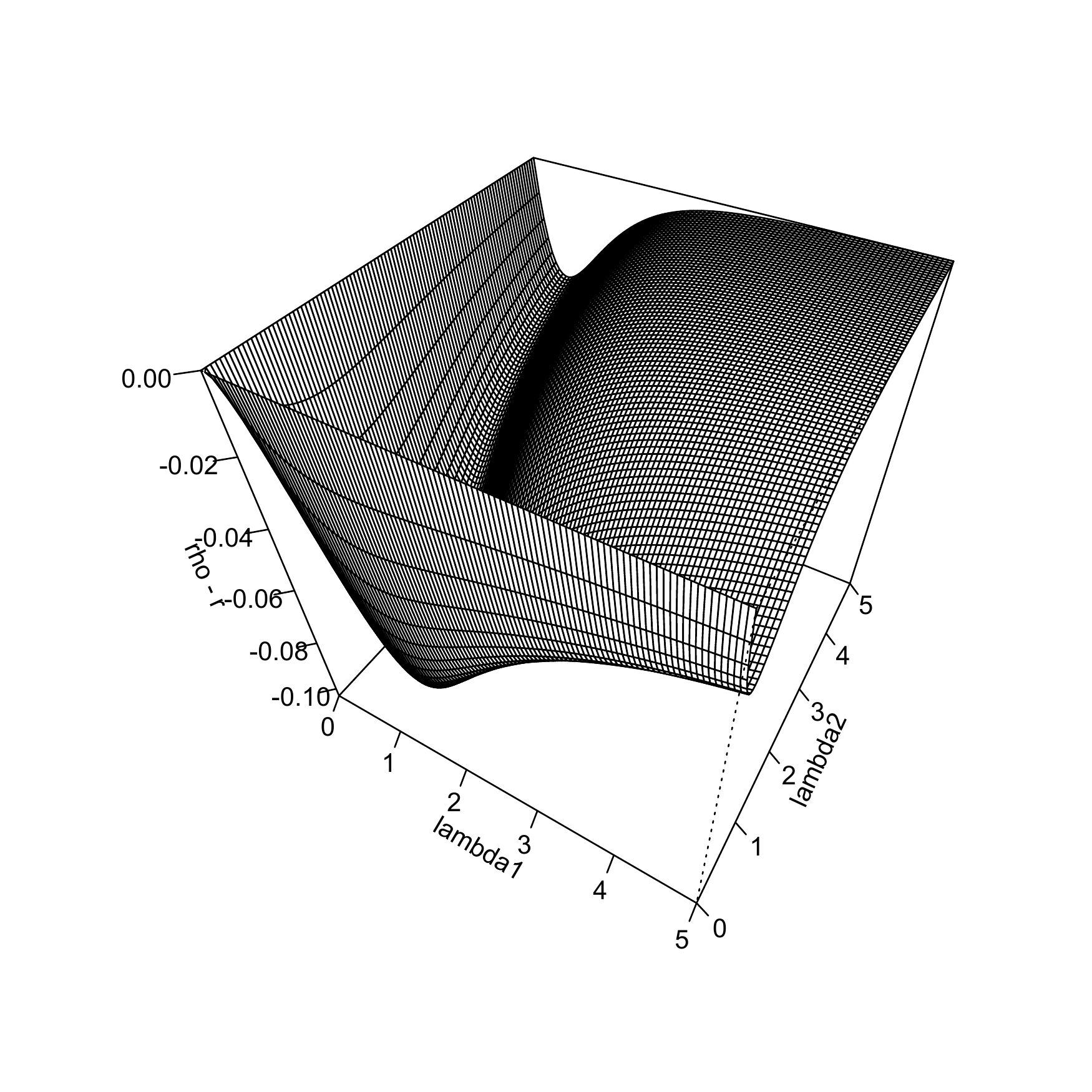}}  \includegraphics[width=6cm]{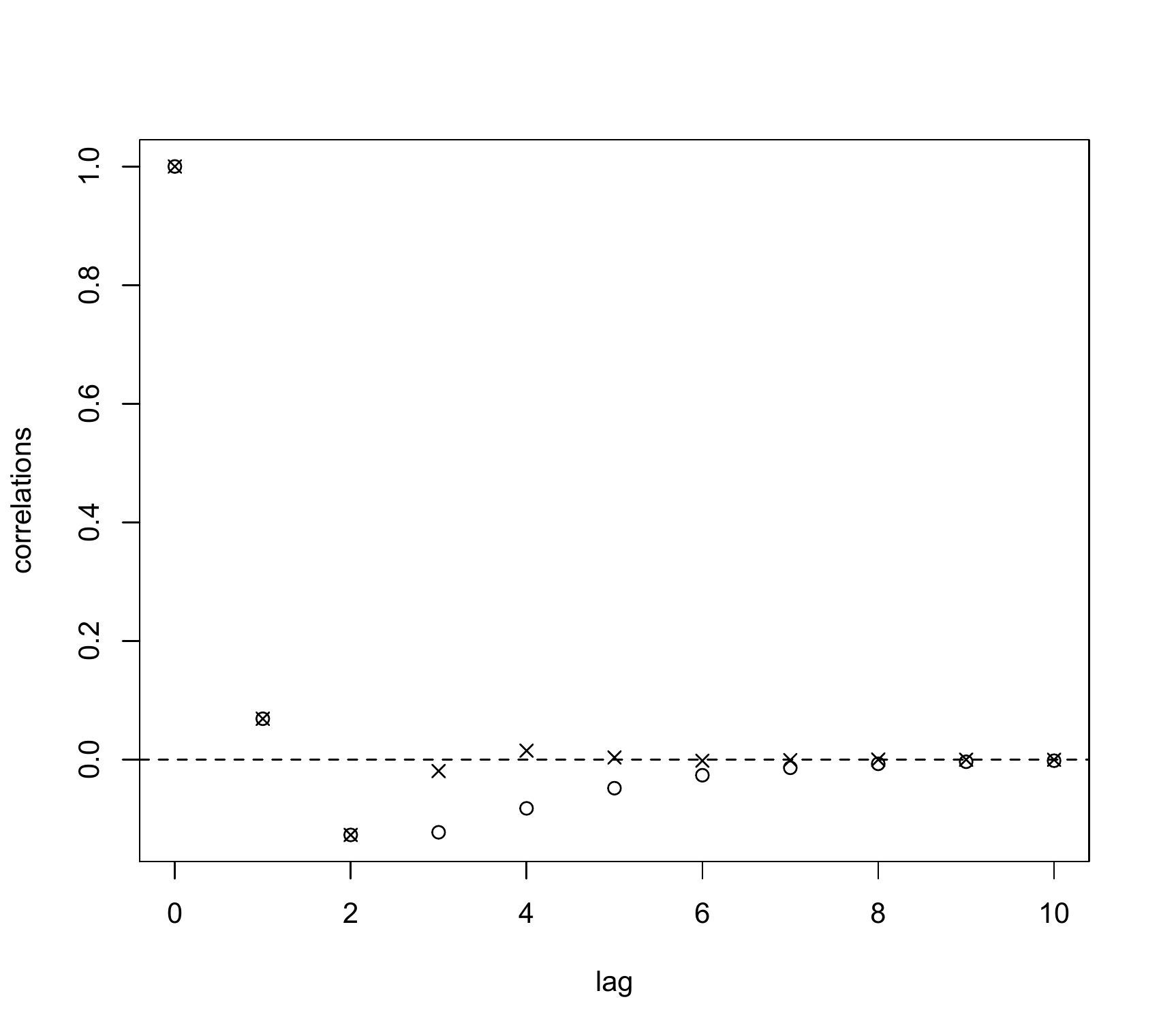}

\vspace{-1cm}

\centerline{(a)\hspace{6cm}(b)}
\vspace{-.2cm}

\caption{(a) Differences of the third order correlations of OU(2) processes with real parameters $\lambda_1,\lambda_2$ and the AR(2) processes with the same first two correlations. (b) The first ten correlations of  two particular OU(2) ($\circ$) and AR(2) ($\times$) processes with $r_1=\rho_1$ and $r_2=\rho_2$. }\label{figdiffer}

\end{figure}

The substitution  $\lambda_h/\tau$ for $\lambda_h$ 
shows that the same family of correlations is  obtained for any value of $\tau$; hence, we can fix $\tau=1$ without loss of generality.
In particular, for $\rho_1=r_1$ and $\rho_2=r_2$, we compute 
 $\rho_3$ and $r_3$ as functions of  $\lambda_1$ and $\lambda_2$ using (\ref{lar3}) and (\ref{larho}).

 The plot of  $\rho_3-r_3$ in Figure \ref{figdiffer} (a) shows that these results differ. In particular, Figure \ref{figdiffer} (b) shows that the correlations of a certain OU(2) ($\circ$) and the AR(2) ($\times$) process with $r_1=\rho_1$ and $r_2=\rho_2$ are not identical. These processes have been chosen to maximise $|\rho_3-r_3|$. The parameters are $\bm\kappa=(0.84,0.84)$ and $r_3-\rho_3=0.1032608$.

\section{Estimation of the parameters of OU($p$)}\label{proce}

\subsection{Reparameterisation by means of real parameters}

Since we wish to consider real processes $x$ and the process itself and its covariance $\gamma(t)$ depend only on the unordered set of the components of 
$\bm\kappa$, 
we shall reparameterise the process by means of the real vector $\bm\phi=(\phi_1,\dots,\phi_p)$ given by the polynomial identity
\begin{equation}\label{polinv}
g(z)=\prod_{j=1}^p(1+\kappa_jz)=1-\sum_{j=1}^p\phi_jz^j.
\end{equation}
The resulting process is real, because of Corollary \ref{rea}.

\subsection{Maximum likelihood estimation (MLE)}\label{maxver}

We shall assume that the process $\mu+x$ is observed at times $0,\tau,2\tau,\dots,n\tau$. By choosing $\tau$ the time unit of measure, we assume without loss of generality that our observations are $\{\mu+x(i): i=0,1,\dots,n\}$. 

The likelihood $L$ of the vector $\Delta\bm x=(x(1)-x(0),x(2)-x(1),\dots,x(n)-x(n-1))$
 is given by
$$\log L(\bm x;\bm\phi,\sigma)=-{\textstyle\frac{n}{2}}\log(2\pi)-{\textstyle\frac{1}{2}}\log(\det(V(\bm\phi,\sigma))-{\textstyle\frac{1}{2}}\Delta\bm x\tras (V(\bm\phi,\sigma))^{-1}\Delta\bm x$$ 
with $V(\bm\phi,\sigma)$ equal to the $n\times n$ matrix with components 
$$V_{h,i}=2\gamma(|h-i|)-\gamma(|h-i|+1)-\gamma(|h-i|-1)$$
 that reduce to $2(\gamma(0)-\gamma(1))$ at the diagonal $h=i$.

From these elements, a numerical optimisation leads to obtain the maximum likelihood estimators $\hat{\bm\phi}$ of $\bm\phi$ and $\hat\sigma^2$ of $\sigma^2$. If required, the estimations  $\hat{\bm\kappa}$ follow by solving the analogue of the polynomial equation  (\ref{polinv}) written in terms of the estimators:
$$\prod_{j=1}^p(1+\hat\kappa_jz)=1-\sum_{j=1}^p\hat\phi_jz^j.$$

The optimisation for large $n$ and the solution of the algebraic equation for large $p$ require a considerable computation effort, but there are efficient programs to perform both operations, as {\sf optim} and {\sf polyroot} in R (\cite{R}). 

An alternative when the process is assumed to be centred ($\mu=0$) is to maximise the log-likelihood of $\bm x=(x(0), x(1),\dots,x(n))^{\tr}$
$$\log L(\bm x;\bm\phi,\sigma)=-{\textstyle\frac{n}{2}}\log(2\pi)-{\textstyle\frac{1}{2}}\log(\det(\Gamma(\bm\phi,\sigma))-{\textstyle\frac{1}{2}}\bm x\tras (\Gamma(\bm\phi,\sigma))^{-1}\bm x$$ 
where $\Gamma$ has components 
$\Gamma_{h,i}=\gamma(|h-i|)$ ($h,i=0,1,\dots,n$).

The optimisation procedures require an initial guess about the value of the parameter to be estimated. The estimators obtained by  {\em matching correlations} described in the next section can be used for that purpose.

\subsection{Matching correlations estimation (MCE)}\label{match}

From the closed formula for the covariance $\gamma$ (eq. (\ref{cf})) and the relationship between $\bm \kappa$ and $\bm\phi$ (eq. (\ref{polinv})), we have  a mapping 
$(\bm\phi,\sigma^2)\mapsto\gamma(t)$, for each $t$.
Since $\bm\rho^{(T)}:=(\rho(1),\rho(2),\dots,\rho(T))\tras$ $=(\gamma(1),\gamma(2),\dots,\gamma(T))\tras/\gamma(0)$ does not depend on $\sigma^2$, these equations determine a map 
${\cal C}:(\bm\phi,T)\mapsto \bm\rho^{(T)}={\cal C}(\bm\phi, T)$ for each $T$.
After choosing  a value of $T$ and obtaining an estimate 
$\bm\rho_e^{(T)}$ of $\bm\rho^{(T)}$ based on the empirical covariances of $x$, 
we propose as a first estimate of $\bm\phi$, the vector  $\check{\bm\phi}_T$ such that all the components of the corresponding $\bm\kappa$ have positive real parts, and such that the euclidean norm $\|\bm\rho_e^{(T)}-{\cal C}(\check{\bm\phi}_T,T)\|$ reaches its minimum.
The procedure resembles the estimation by the {\em method of moments}.
The components of $\bm\rho_e^{(T)}$for the series $(x_i)_{i=1,2,\dots,n}$ are computed as $\rho_{e,h}=\gamma_{e,h}/\gamma_{e,0}$, $\gamma_{e,h}=\frac{1}{n}\sum_{i=1}^{n-h} x_ix_{i+h}$.

\subsection{Some simulations}

We have simulated the series $x(i),i=0,1,2,\dots,n$ obtained from an OU process $x$ for $n=300$ and three different values of the parameters and computed the MC and ML estimators $\check{\bm\phi}_T$, and $\hat{\bm\phi}$. The value of $T$ for the MC estimation has been arbitrarily set equal to the integral part of $0.9\times n$, but the graphs of $\check{\bm\phi}_T$ for several values of T show in each case that after $T$ exceeds a moderate threshold, the estimates remain practically constant. One of such graphs is included below (see Figure \ref{oupppej1}).

The simulations show that the correlations of the series with the estimated parameters are fairly adapted to each other and to the empirical covariances. The departure from the theoretical covariances of $x$ can be ascribed to the simulation intrinsic randomness.

Our first two examples
describe OU(3) processes with arbitrarily (and randomly) chosen parameters and the third one imitates the behaviour of Series A that appears in \S\ref{examp}.

\begin{example}\label{ejem1} 
A series $(x_i)_{i=0,1,\dots,n}$ of $n=300$ observations of the OU$_{\sbm\kappa}$ process $x$ ($p=3$, $\bm\kappa=(0.9, 0.2+0.4\imath 0.2-0.4\imath)$, $\sigma^2=1$) was simulated, and the parameters $\bm\phi=(-1.30$, $-0.56$, $-0.18)$ and $\sigma^2=1$ were estimated by means of $\check{\bm\phi}_T$ $=(-1.9245,$ $ -0.6678, $ $-0.3221)$, $T$ $=270$,  $\hat{\bm\phi}$ $=(-1.3546$, $-0.6707$, $ -0.2355)$ and $\hat\sigma^2=0.8958$.
The corresponding estimators for $\bm\kappa$ are $\check{\bm\kappa}$ $=(1.6368,$ $ 0.1439$ $+0.4196\imath, $ $0.14389$ $-0.4196\imath)$ and $\hat{\bm\kappa}$ $=(0.9001$, $0.2273+0.4582\imath$, $0.2273-0.4582\imath)$.

\begin{figure}[!h]
\centerline{\includegraphics[width=8cm]{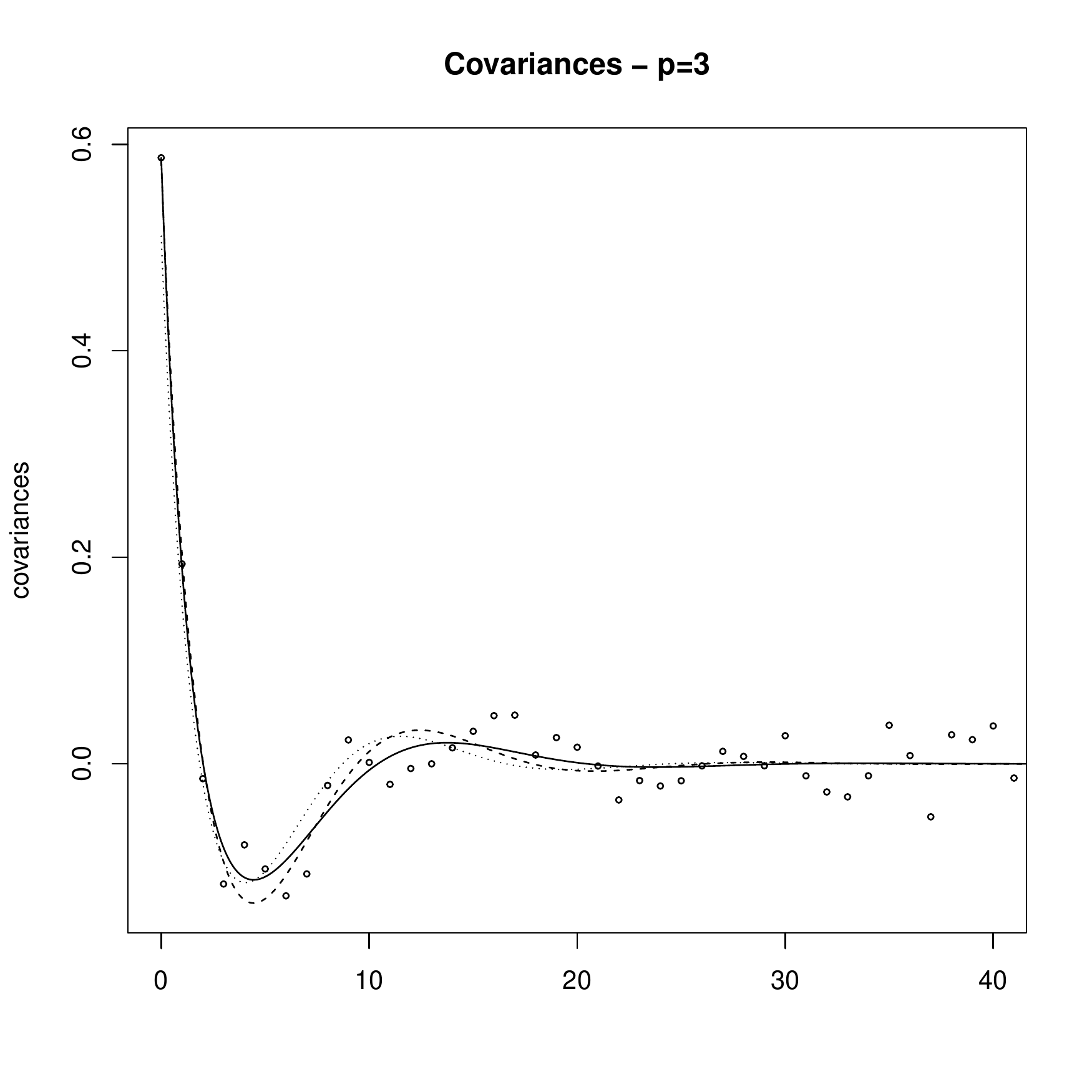}}
\caption{Empirical covariances ($\circ$) and covariances of the MC (---) and ML (- - -) fitted OU models, for $p=3$ corresponding to Example 1. The covariances of OU$_{\sbm\kappa}$ are indicated with a dotted line.}\label{oupej1}
\end{figure}

Figure \ref{oupej1} describes the theoretical, empirical and estimated covariances of $x$ under the assumption $p=3$, that is, the actual order of $x$. The results obtained when the estimation is performed for $p=2$ and $p=4$ are shown in Figure \ref{ouppej1}. Finally, Figure \ref{oupppej1} shows that the MC estimates of $\bm\phi$ become stable for $T$ moderately large, and close to the already indicated estimations for $T=270$ (the horizontal lines). 

\begin{figure}[!h]
\centerline{\includegraphics[width=6cm]{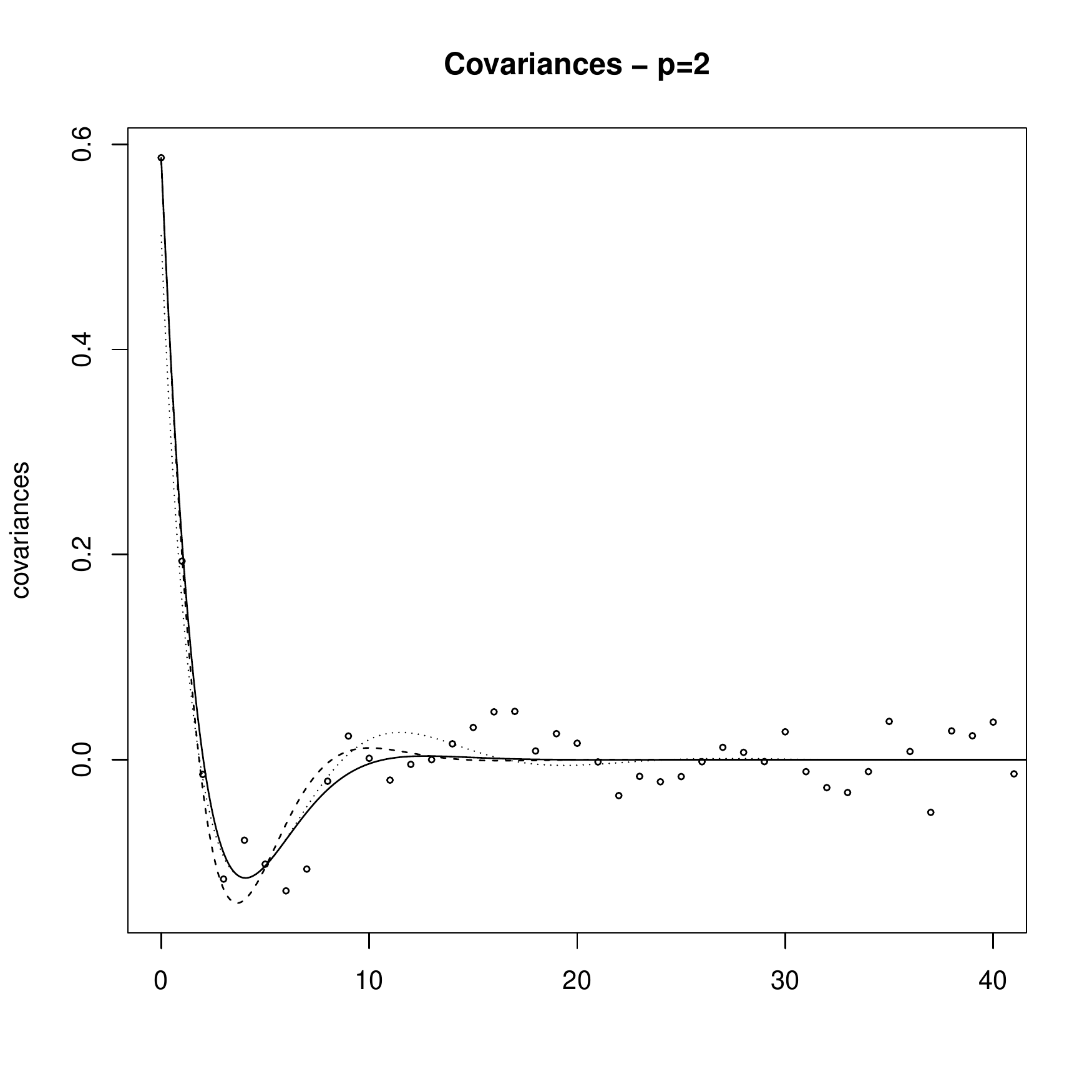}\includegraphics[width=6cm]{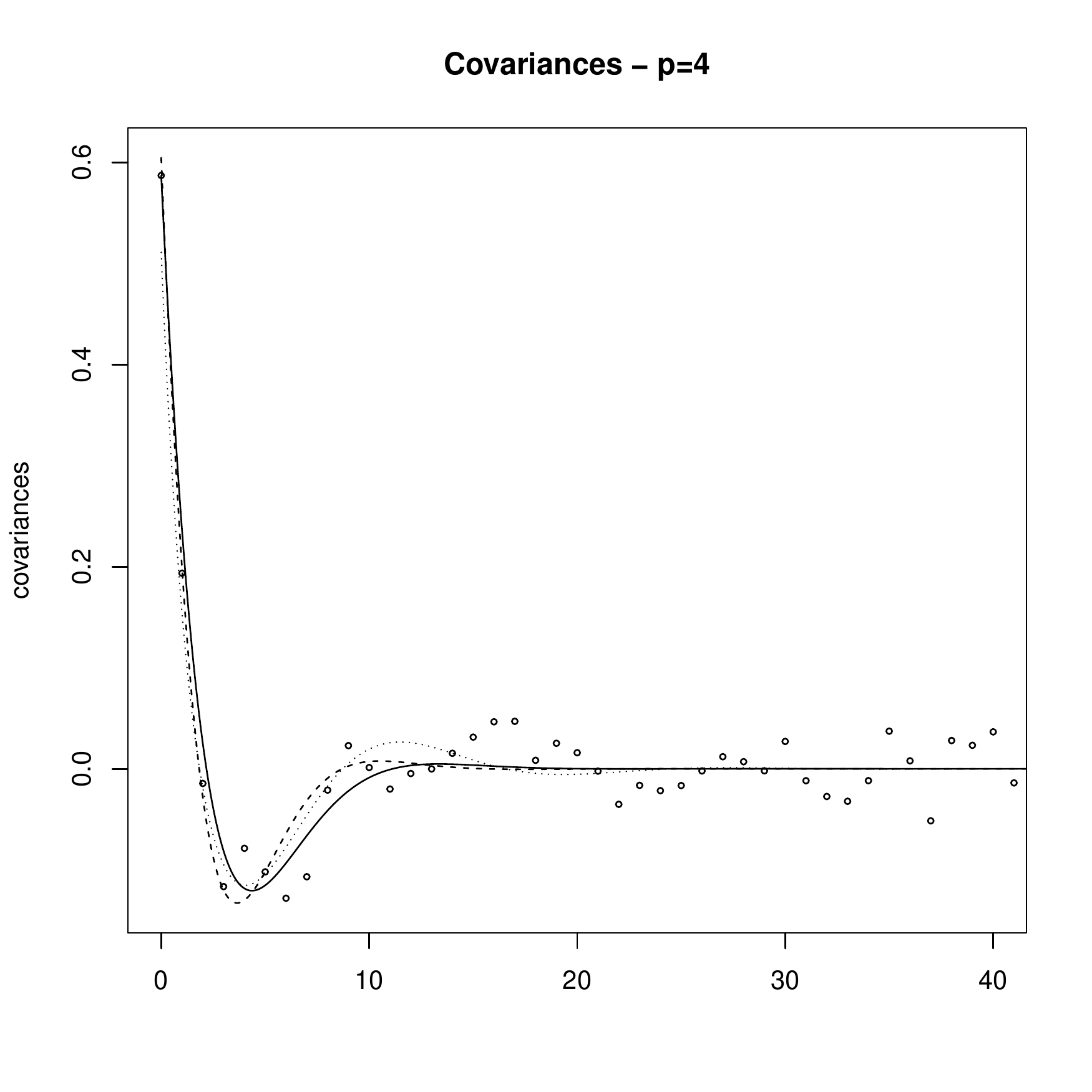}}

\vspace{-4mm}

\caption{Empirical covariances ($\circ$) and covariances of the MC (---) and ML (- - -) fitted OU models, for $p=2,4$ corresponding to Example 1. The covariances of OU$_{\sbm\kappa}$ are indicated with a dotted line.}\label{ouppej1}
\end{figure}

\begin{figure}[!h]
\centerline{\includegraphics[width=6cm]{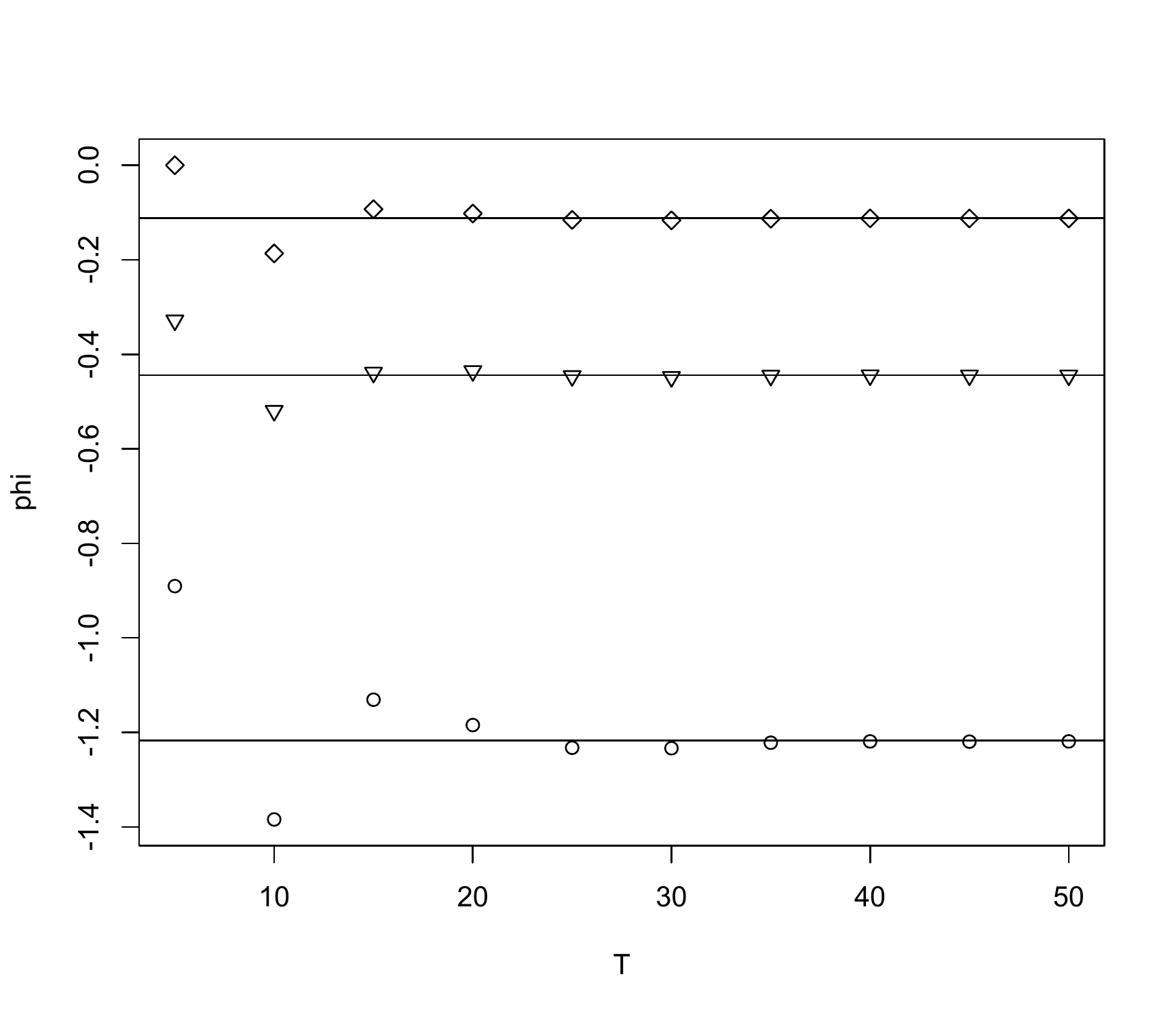}}

\vspace{-4mm}

\caption{The MC estimations $\check\phi_1 (\circ)$, $\check\phi_2 ({\triangledown})$ and $\check\phi_2 (\diamond)$ for different values of $T$, corresponding to Example 1. The horizontal lines indicate the estimations for $T=270$.}\label{oupppej1}
\end{figure}

\end{example} 

\begin{example}\label{ejem2}
The process $x=$ OU$_{(0.04, 0.21, 1.87)}$ is analysed as in Example 1. The resulting estimators are  $\check{\bm\phi}_T=( -2.0611$, $-0.7459$, $-0.0553)$, $T$ $=270$, $\check{\bm\kappa}=(1.6224$, $0.3378$, $0.1009)$, $\hat{\bm\phi}$ $=(-1.8253$, $-0.7340$, $-0.0680)$, $\hat\sigma^2=0.7842$, $\hat{\bm\kappa}$ $=(1.3015$, $0.3897$, $0.1342)$, and the resulting covariances are shown in Figure \ref{oupej2} .

\begin{figure}[!htb]
\centerline{\includegraphics[width=9cm]{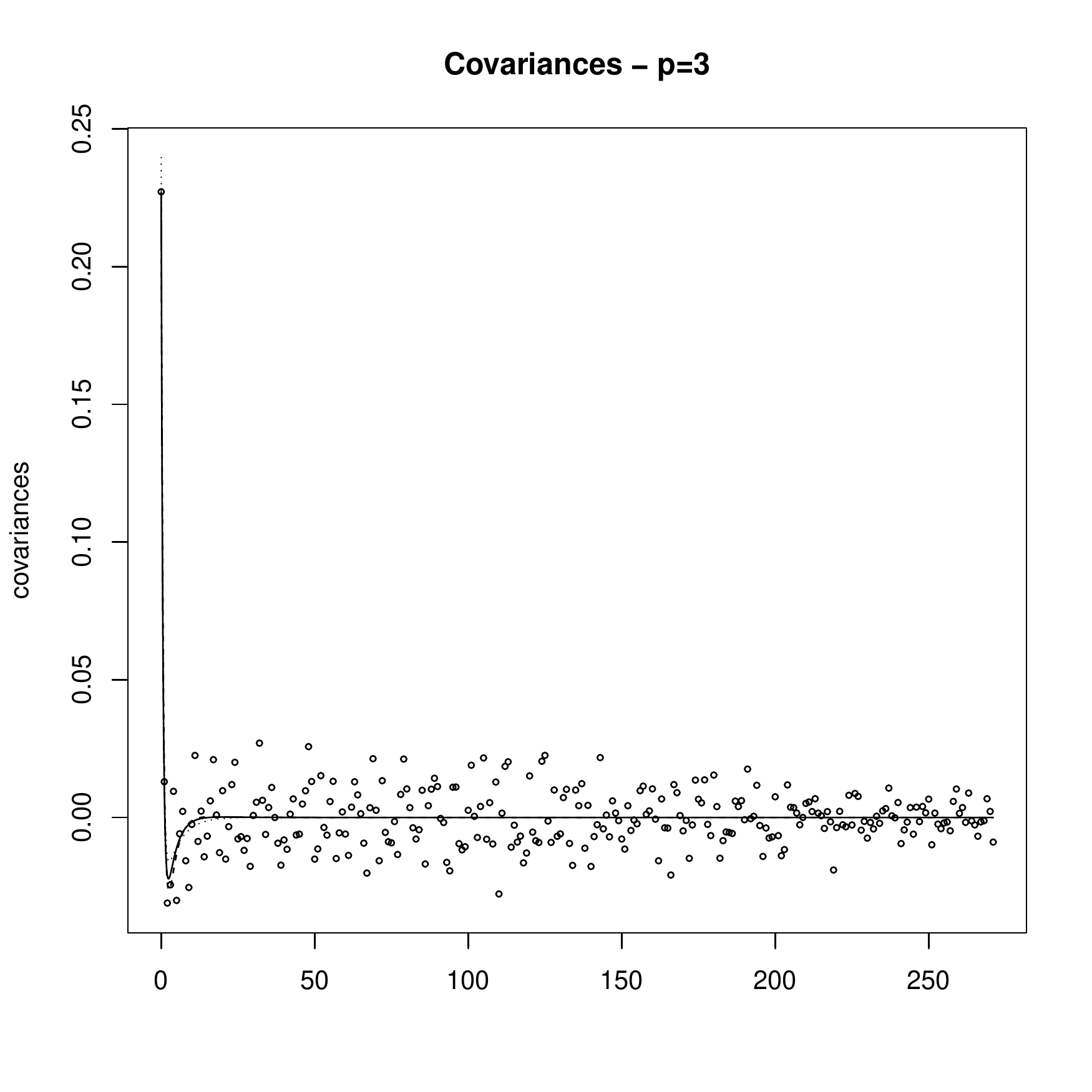}}
\caption{Empirical covariances ($\circ$) and covariances of the MC (---) and ML (- - -) fitted OU models, for $p=3$ corresponding to Example 2. The covariances of OU$_{\sbm\kappa}$ are indicated with a dotted line.}\label{oupej2}
\end{figure}

\end{example}

\begin{example}\label{ejem3}

The parameter $\bm\kappa=(-0.83 -0.0041 -0.0009)$ used in the simulation of the OU process $x$ treated in the present example is approximately equal to the parameter $\hat{\bm\kappa}$ obtained by ML estimation with $p=3$ for Series A in \S\ref{seriesA}. As in previous examples, a graphical presentation of the estimated covariances is given in Figure \ref{oupej3}.

\begin{figure}[!t]

\centerline{\includegraphics[width=4cm]{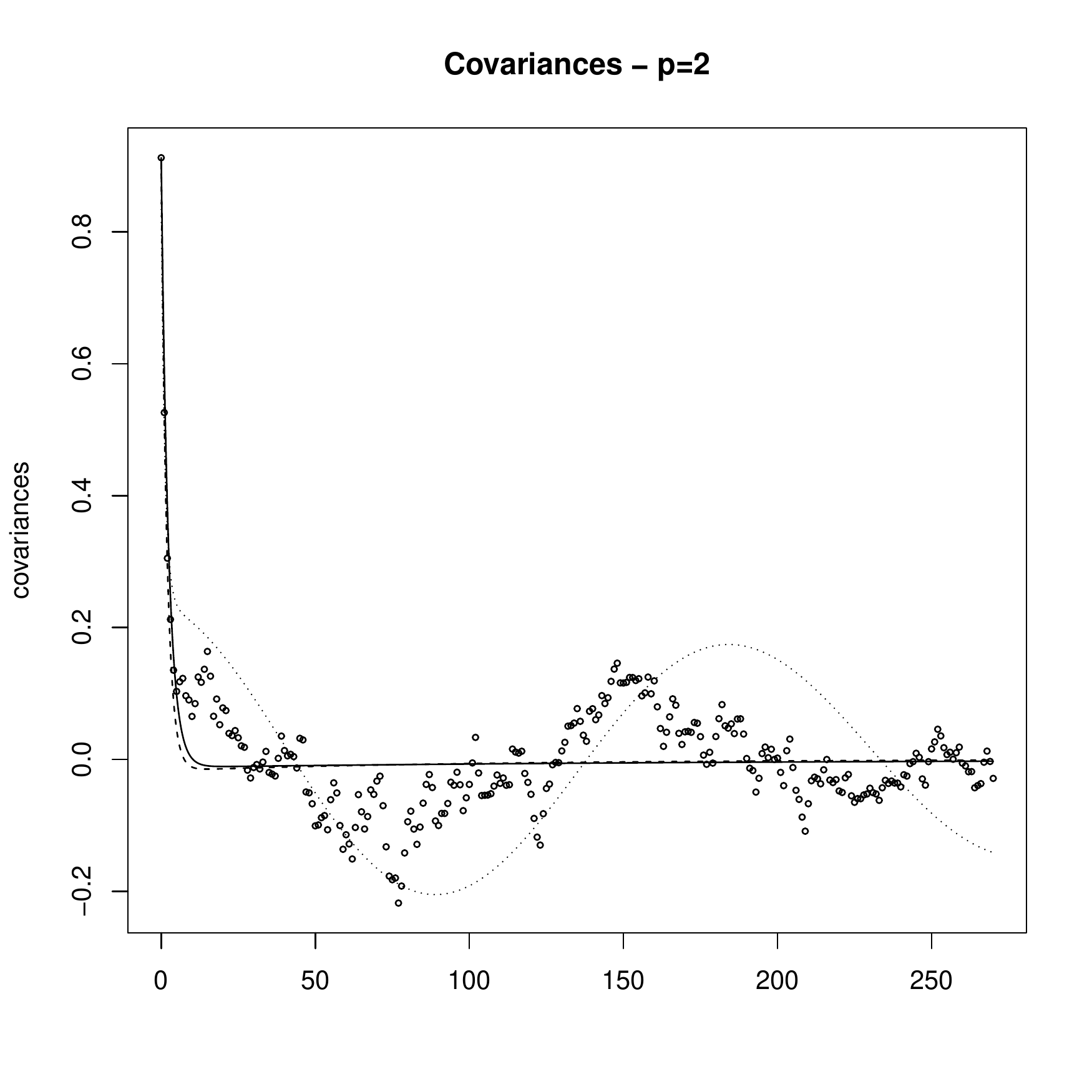}
\includegraphics[width=4cm]{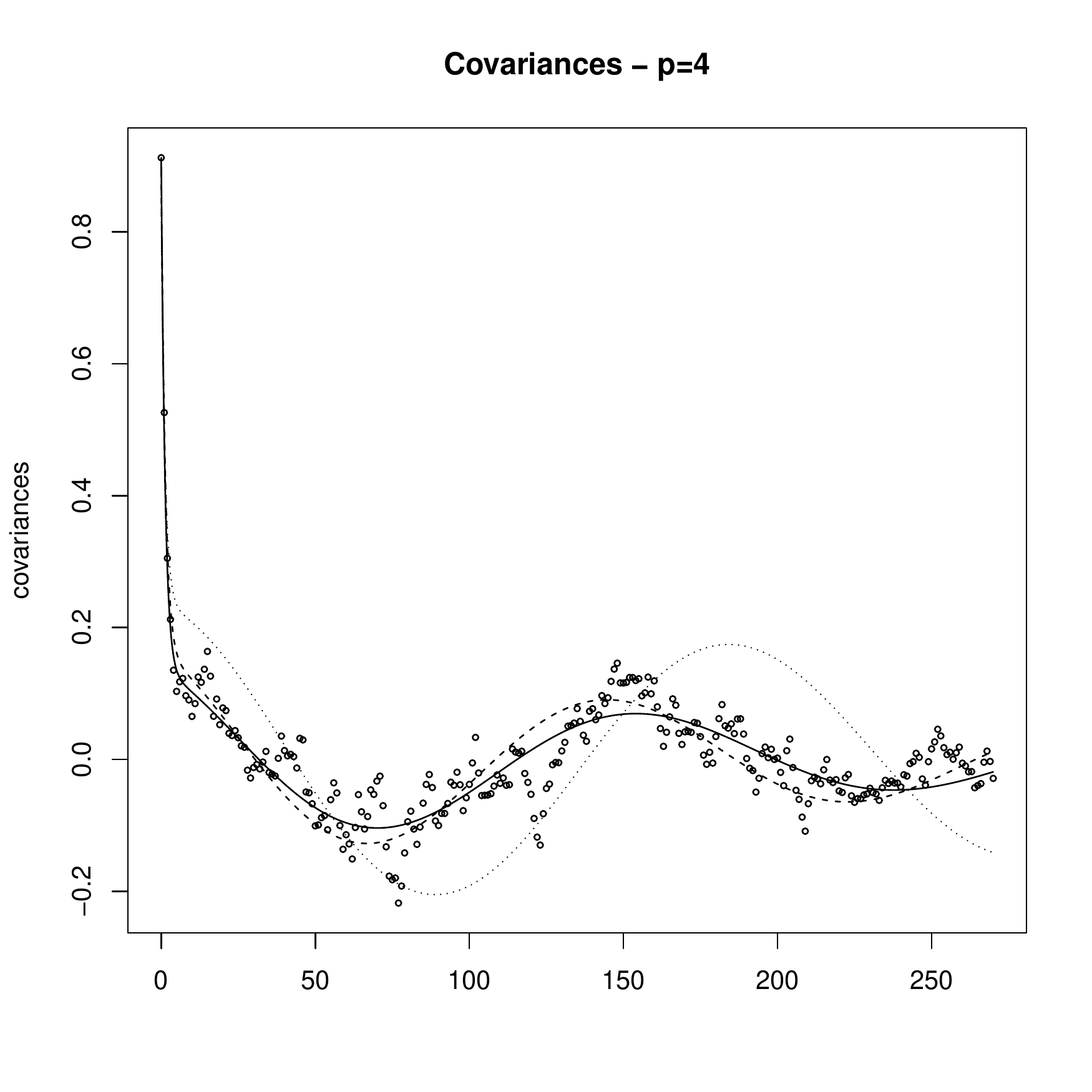}}

\centerline{\includegraphics[width=9cm]{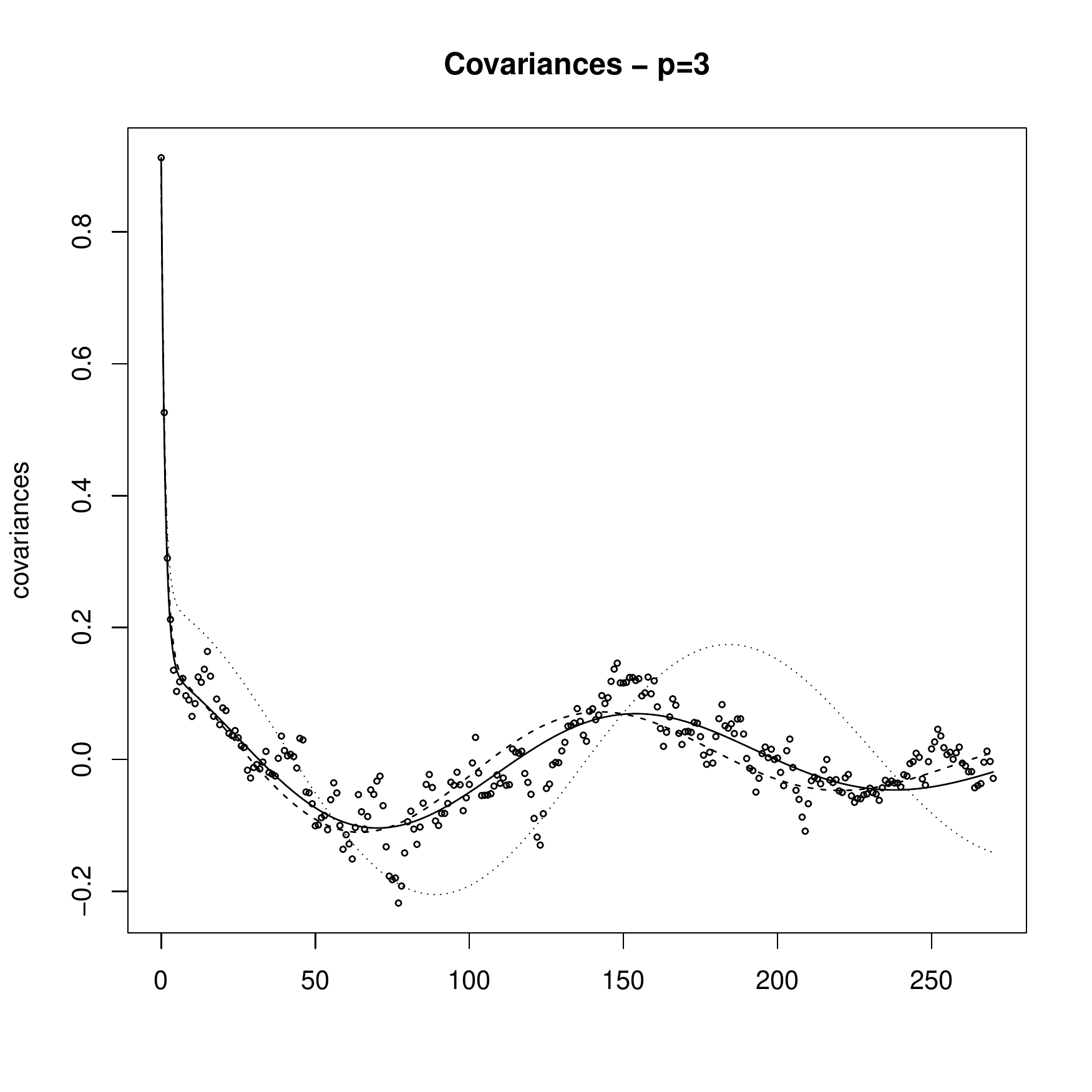}}
\vspace{-8mm}

\caption{Empirical covariances ($\circ$) and covariances of the MC (---) and ML (- - -) fitted OU models, for $p=2, p=4$ and $p=3$, the actual value of the parameter, corresponding to Example 3. The covariances of OU$_{\sbm\kappa}$ are indicated with a dotted line.}\label{oupej3}
\end{figure}
\end{example}

\begin{figure}[!htb]
\centerline{\includegraphics[width=9cm]{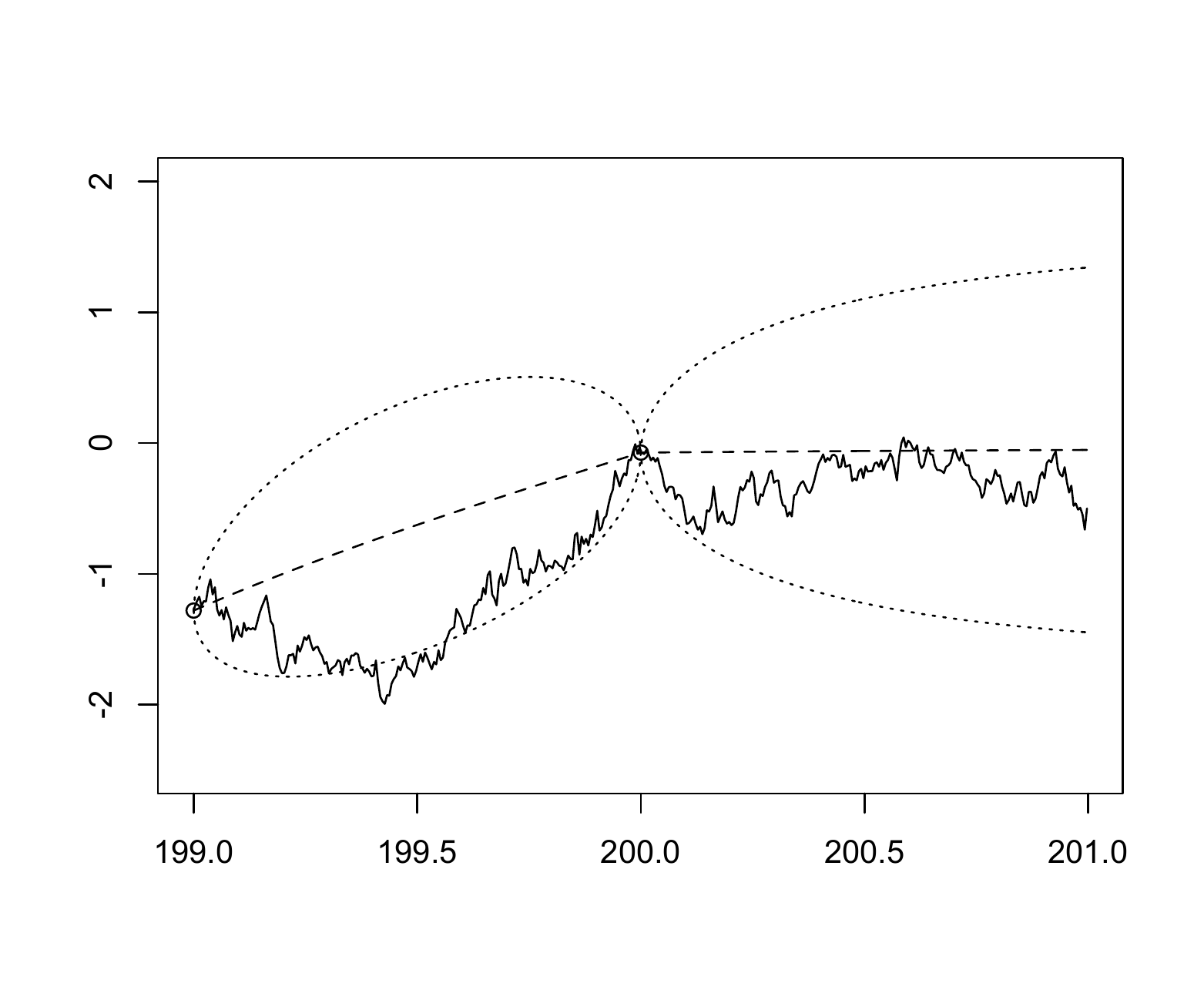}}
\caption{Estimated interpolation and prediction of $x(t)$ for $199<t<200$ and $200<t<201$, respectively (- - -), $2\sigma$ confidence bands based on $(x(i))_{i=0,1,\dots,200}$ ($\cdots$), and a refinement of the simulation of $x(t)$ on $199<t<200$. }\label{interpsimuA}
\end{figure}

The description of the performance of the model is complemented by comparing in Figure \ref{interpsimuA} the simulated values of the process in 400 equally spaced points filling the interval (199,201) with the predicted values for the same interval, based on the OU(3) model and the assumed observed data  $x(0),x(2),x(3),\dots,x(200)$. Also a $2\sigma$ confidence band is included in the graph, in order to describe the precision of the predicted values.

\section{Applications to real data}\label{examp}

In this section we present  experimental results on three real data sets.  We fit  OU$(p)$ processes for small values of $p$ and also some ARMA processes. 
In each case we have observed that  we can find an adequate value of $p$ for which the 
empirical covariances are well approximated by the covariances of the 
adjusted OU$(p)$ model. This is not the case for the AR or ARMA models in all three examples. We present a detailed comparison of both methodologies for the first example.

 The first two data sets are taken from \cite{B&J}, and correspond to equally spaced observations of continuous time processes that might be assumed to be stationary. The third one is a series obtained by choosing one in every 100 terms of a high frequency recording of  oxigen saturation in blood of a newborn child. The data were obtained by a team of researchers of Pereira Rossell Children Hospital in Montevideo, Uruguay, integrated by L. Chiapella,  A. Criado and C. Scavone. Their permission to analyse the data is gratefully acknowledged by the authors.

\subsection{Box, Jenkins and Reinsel ``Series A''}\label{seriesA}

\begin{figure}[!b]
\centerline{\includegraphics[width=6cm]{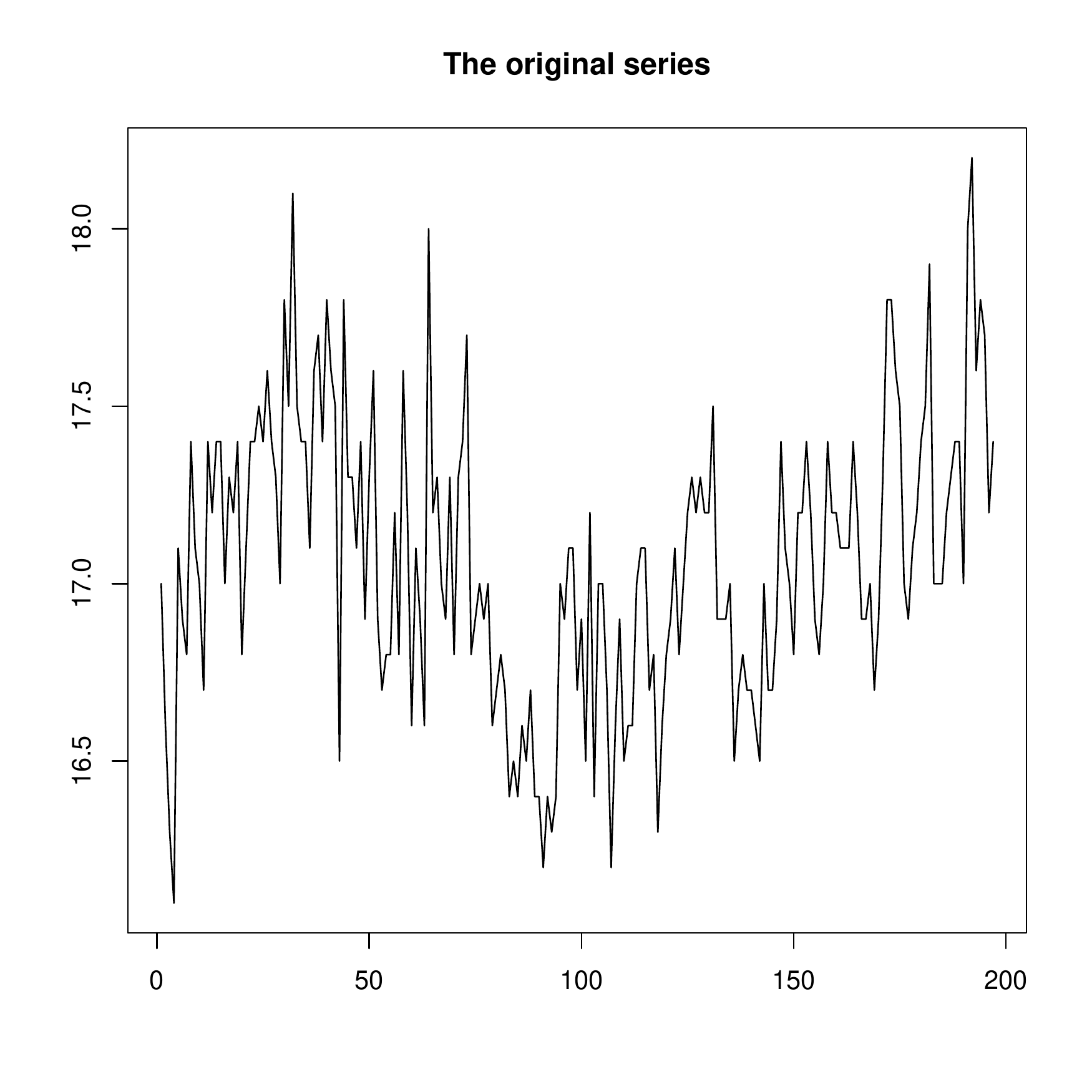}}
\vspace{-8mm}
\caption{Series A}\label{plotsera}
\end{figure}

The Series A is a record of $n=197$ chemical process concentration readings, taken every two hours, introduced with that name and analysed in Chapter 4 of  \cite{B&J} (see also http://rgm2.lab.nig.ac.jp/RGM2/ tfunc.php?rd\underline{ }id=FitAR:SeriesA).  The original data are plotted in Figure \ref{plotsera}.

\begin{figure}[!ht]
\centerline{\includegraphics[width=4cm]{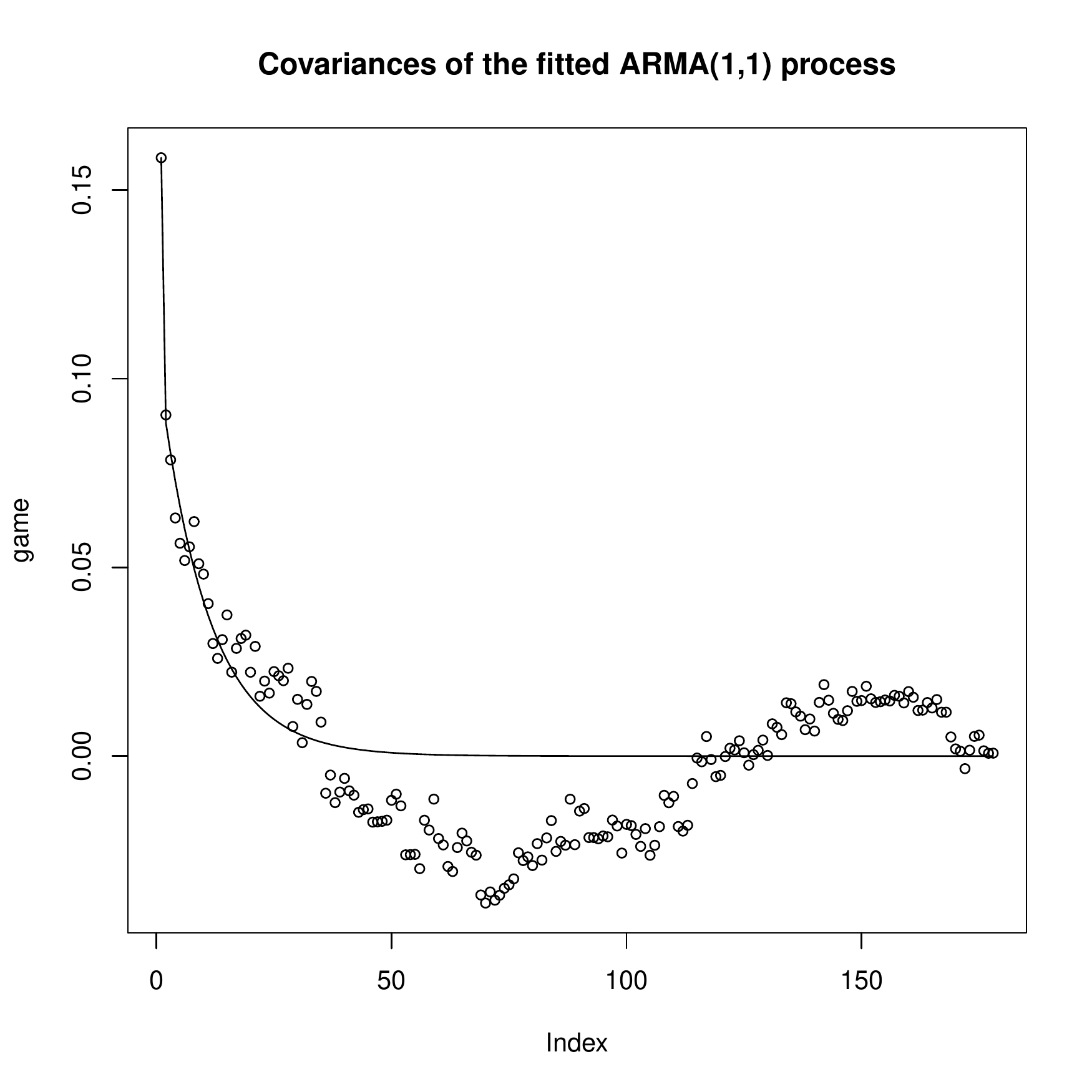}\includegraphics[width=4cm]{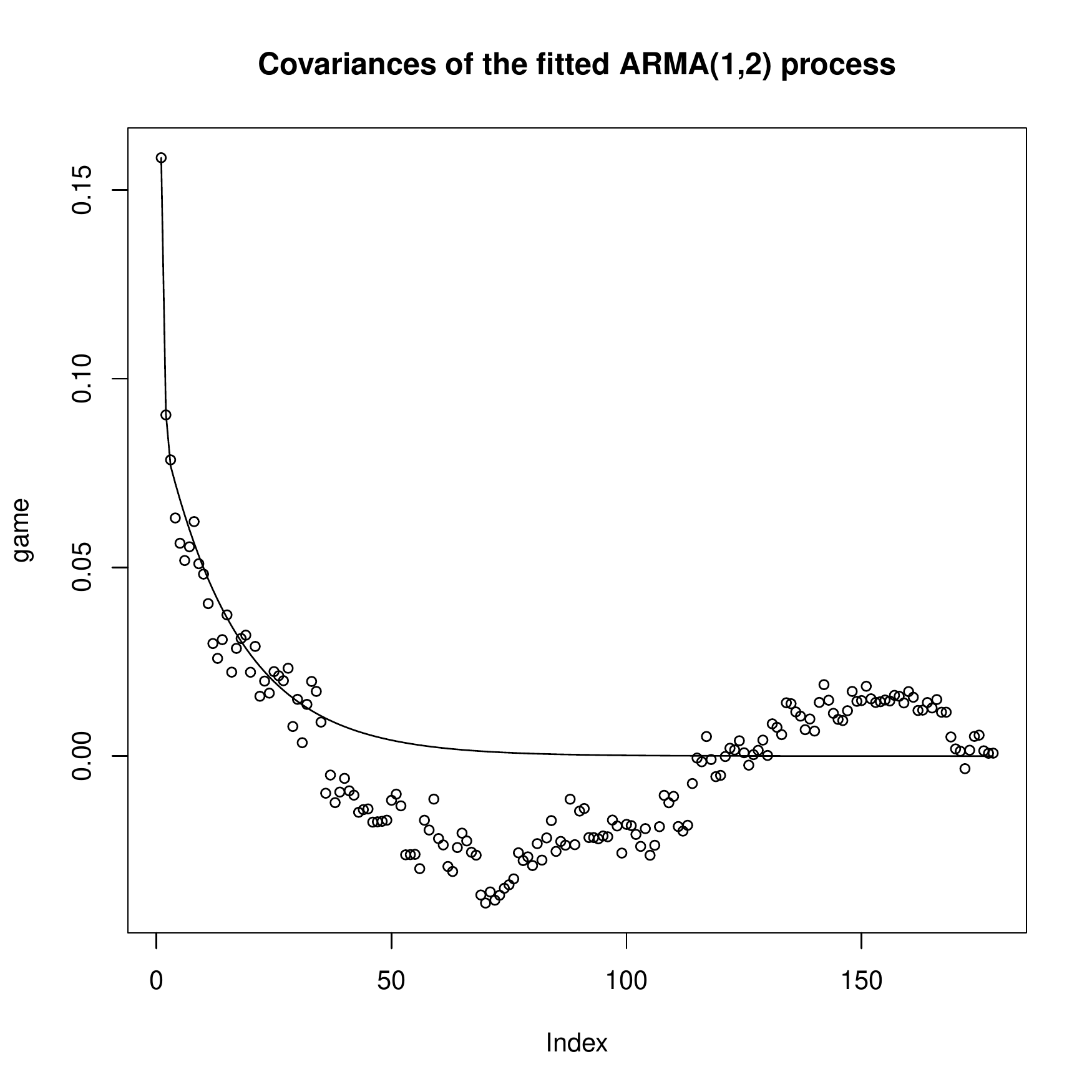}\includegraphics[width=4cm]{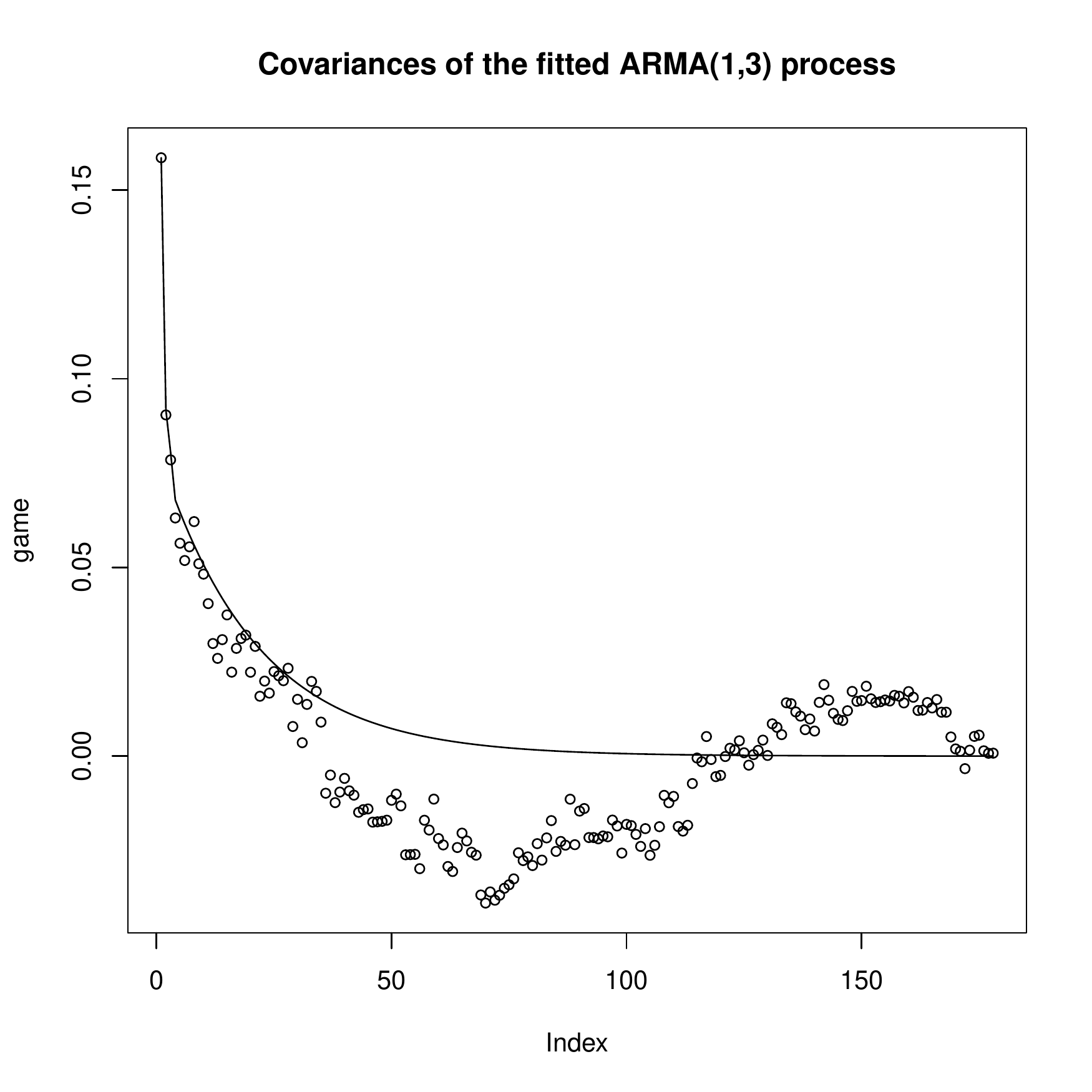}}
\centerline{\includegraphics[width=4cm]{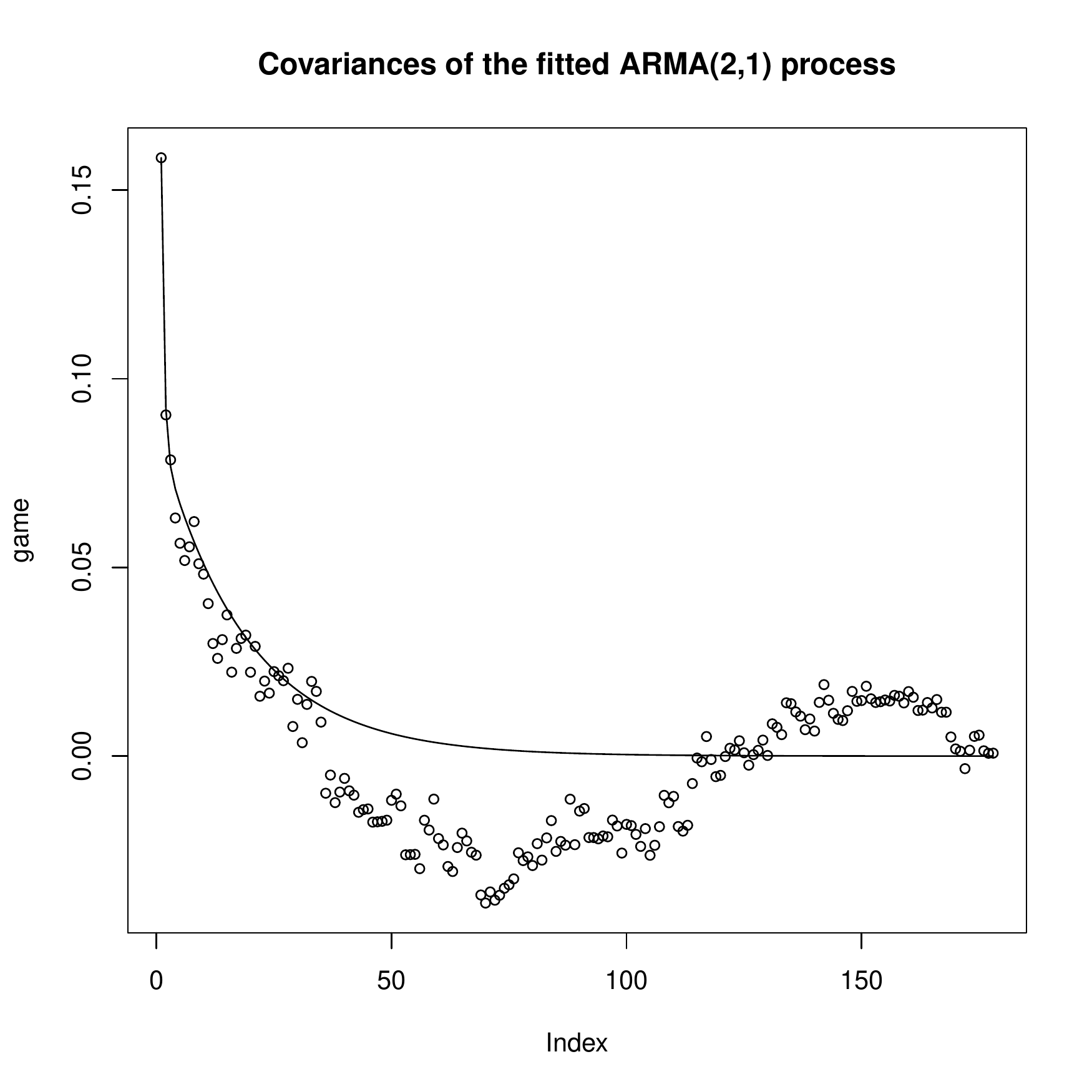}\includegraphics[width=4cm]{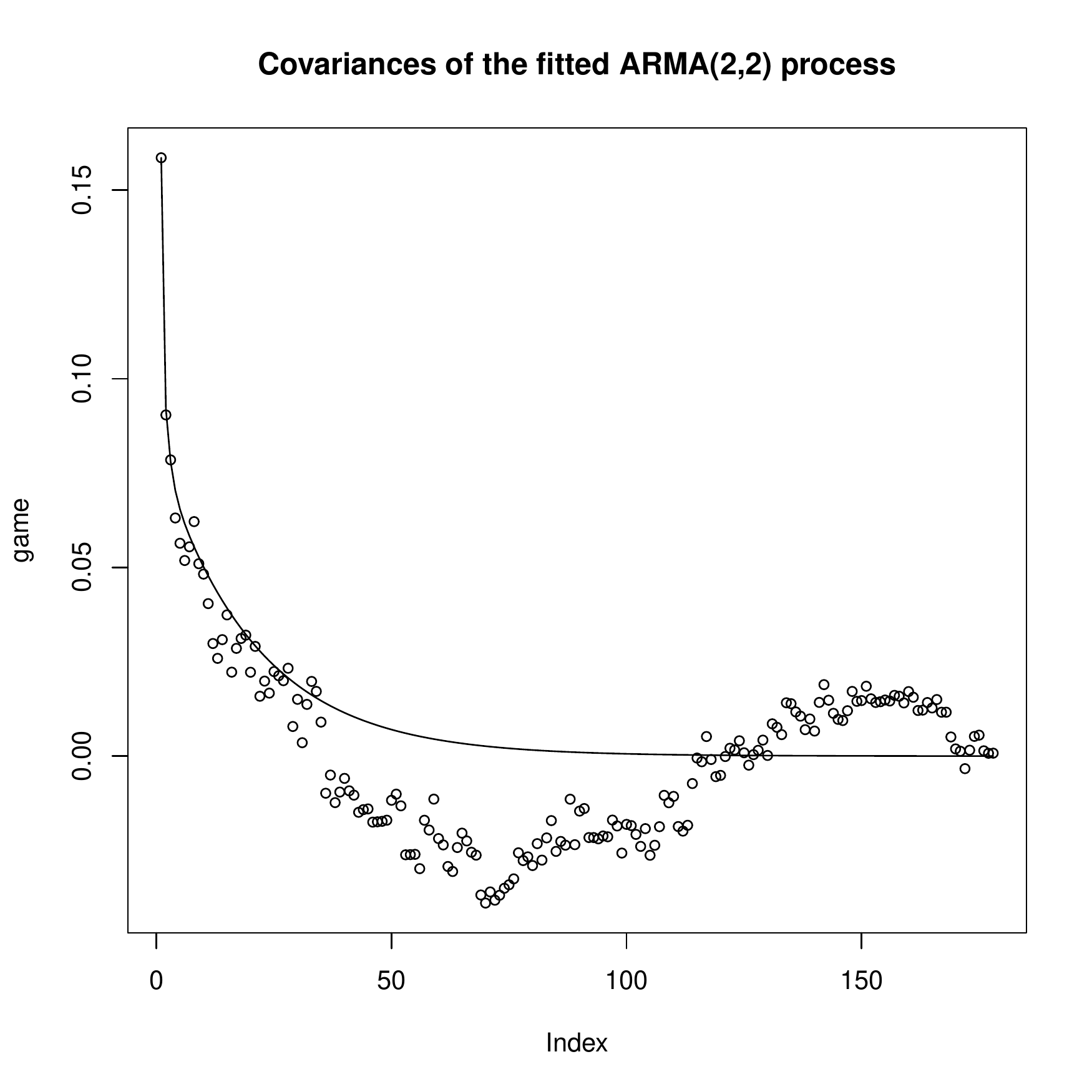}\includegraphics[width=4cm]{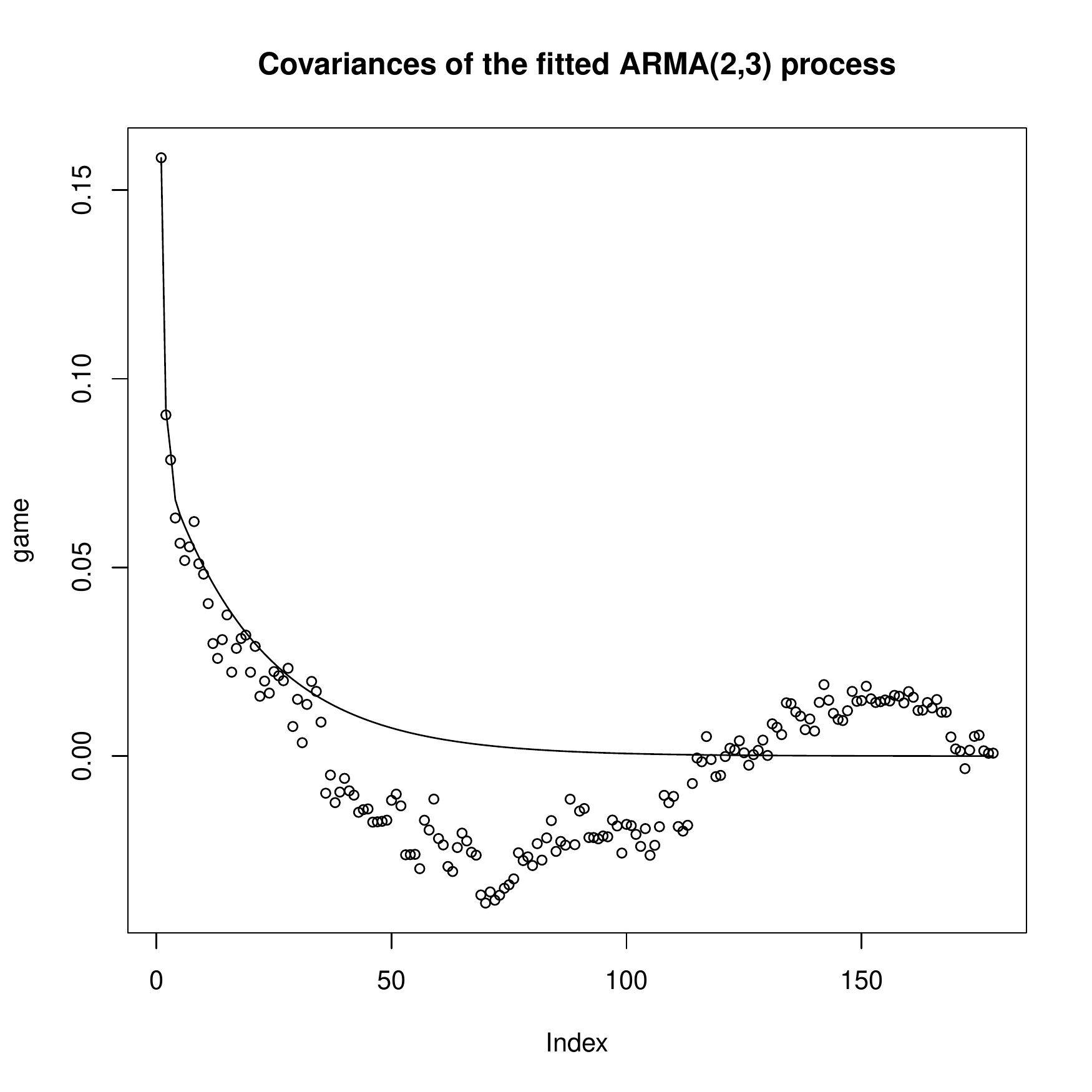}}
\centerline{\includegraphics[width=4cm]{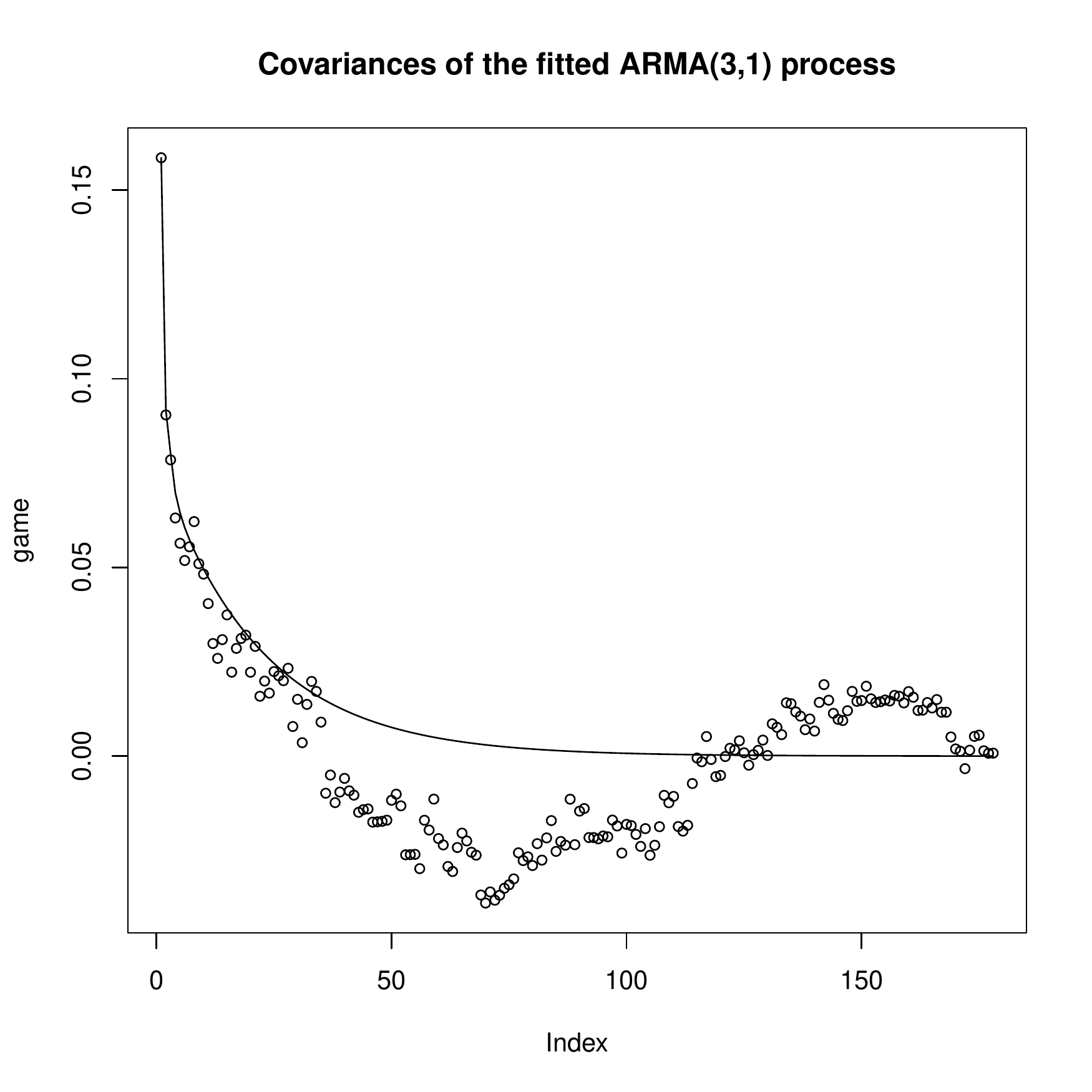}\includegraphics[width=4cm]{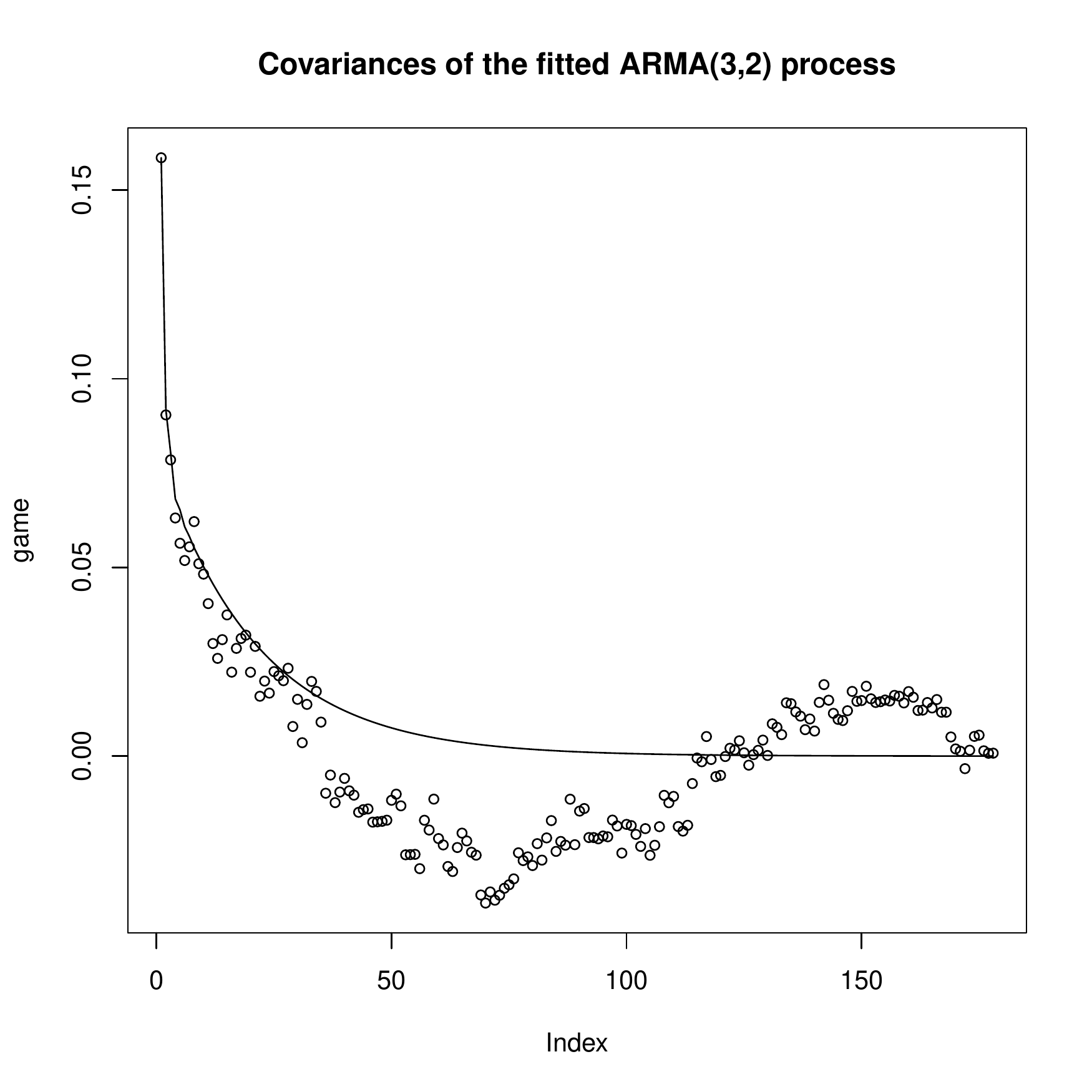}\includegraphics[width=4cm]{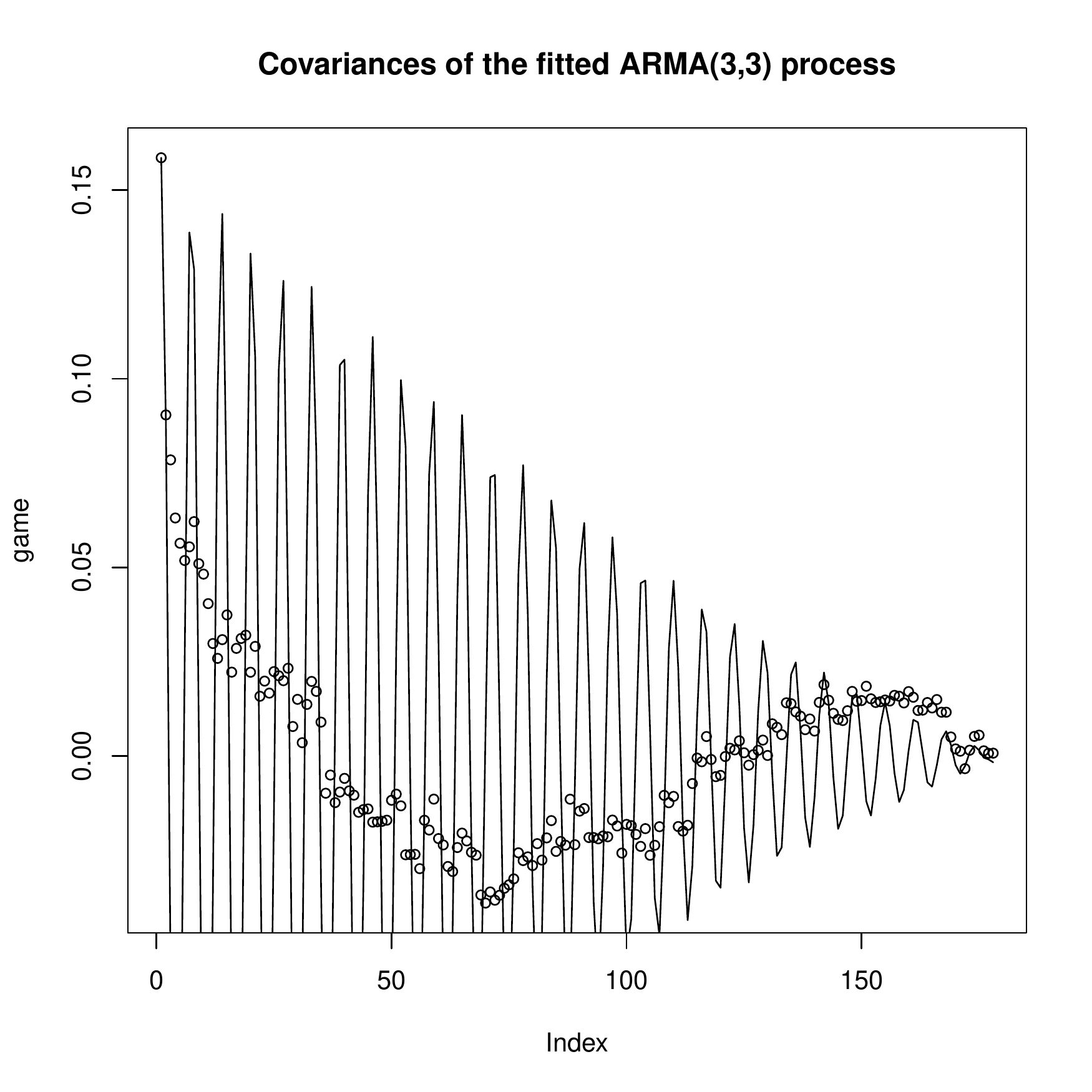}}
\centerline{\includegraphics[width=4cm]{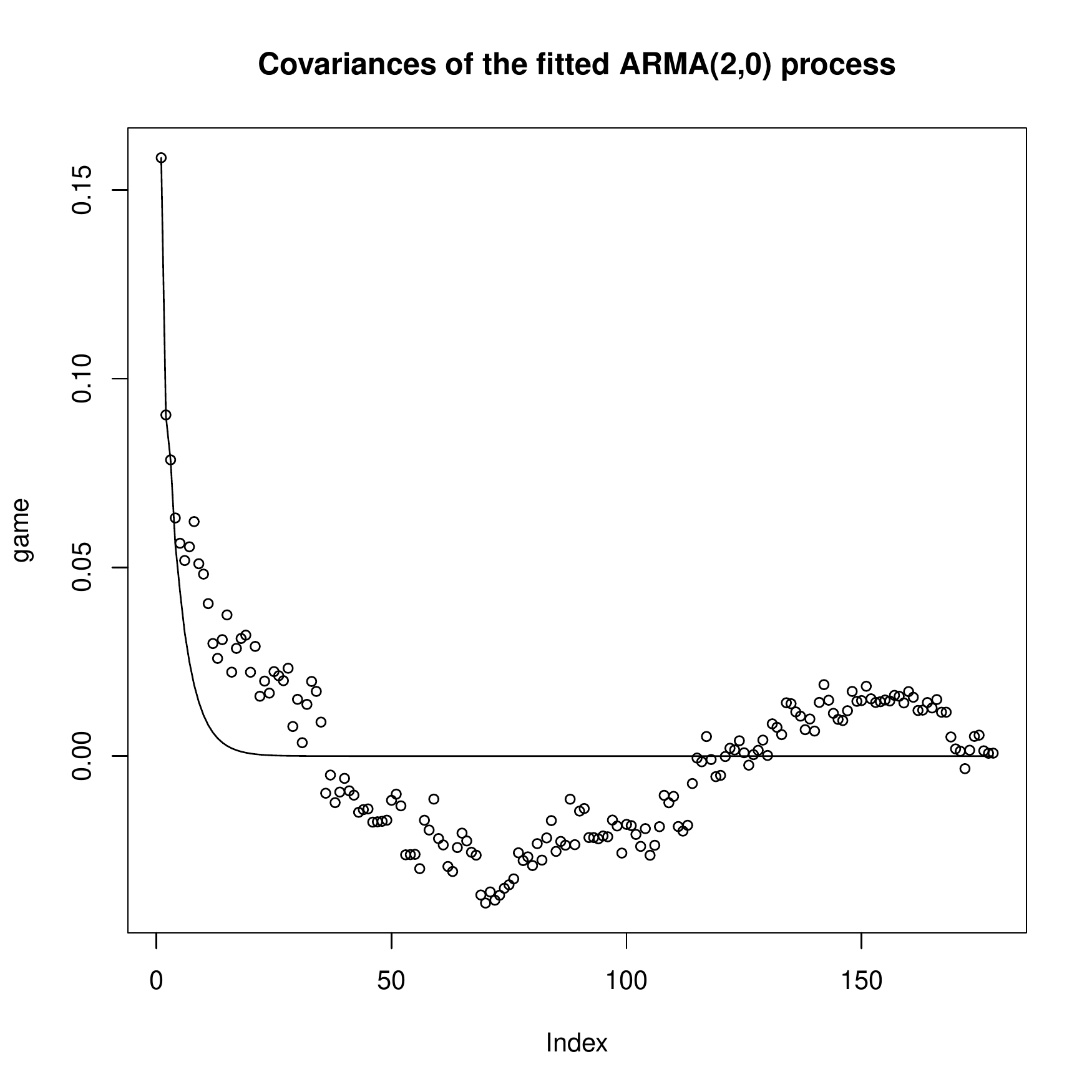}\includegraphics[width=4cm]{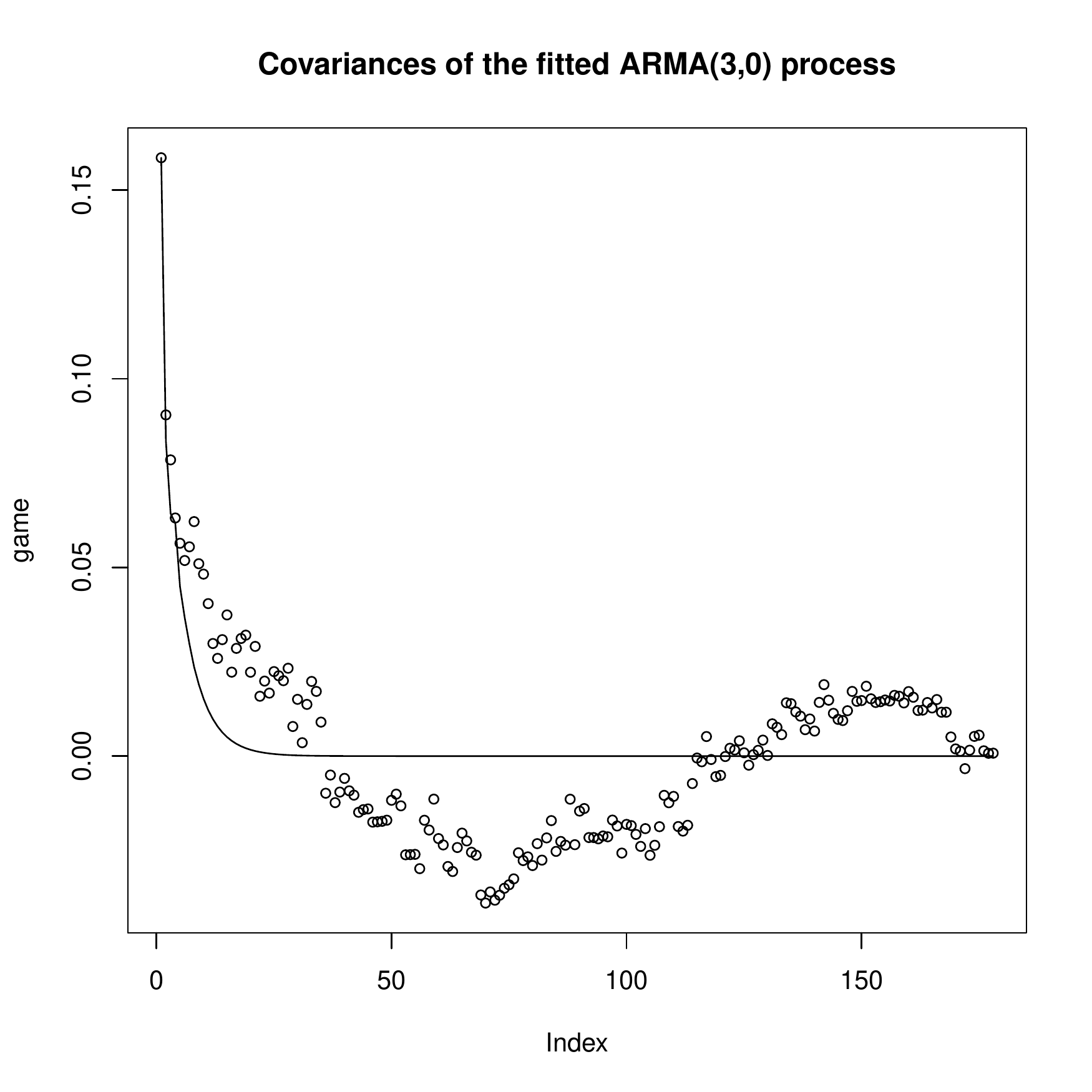}\includegraphics[width=4cm]{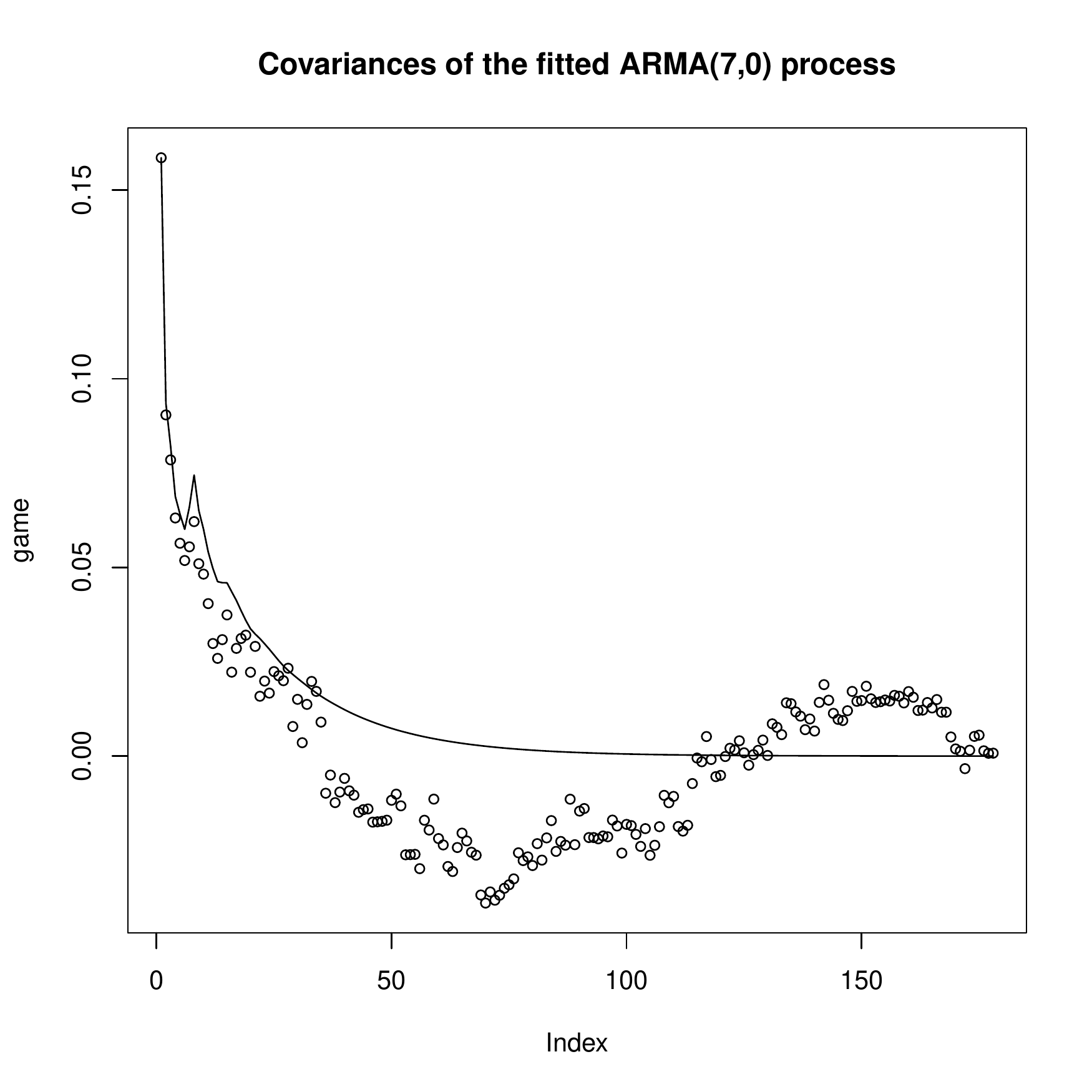}}
\caption{Empirical and ARMA-fitted covariances for Series A}\label{SeriesA/comparma}
\end{figure}

 The diagrams in Figure \ref{SeriesA/comparma}  compare the empirical covariances of the series with the covariances of the estimated 
 ARMA$(p,q)$ process fitted by means of the R function {\sf arima} for several values of $p$ and $q$. In particular, the
ARMA(1,1) is suggested as a model for this data in \cite{B&J}, and subsets of AR(7) are proposed in \cite{cleve} and \cite{McLeod} for the same purpose.

\begin{figure}[!t]
\centerline{\includegraphics[width=4.5cm]{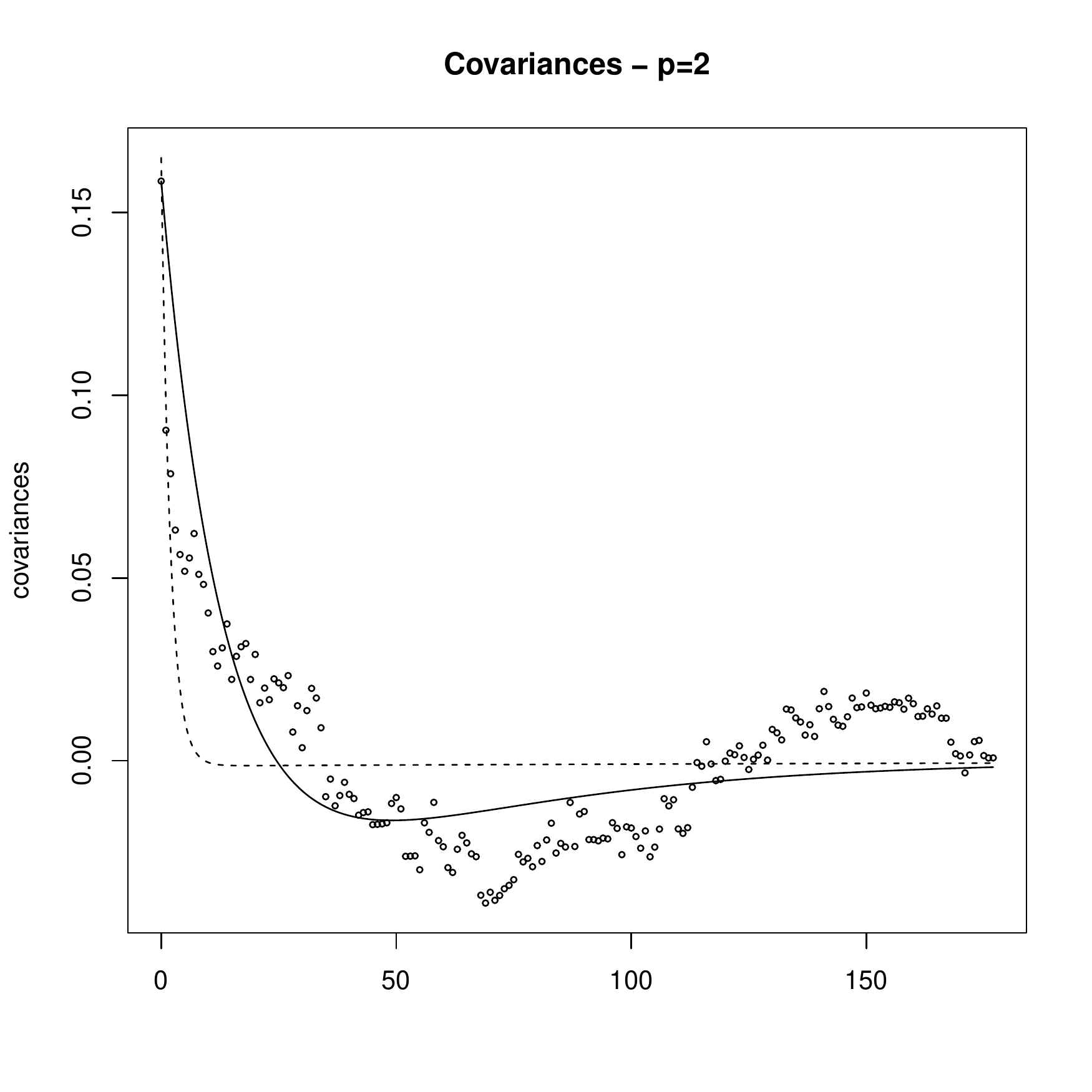}\includegraphics[width=4.5cm]{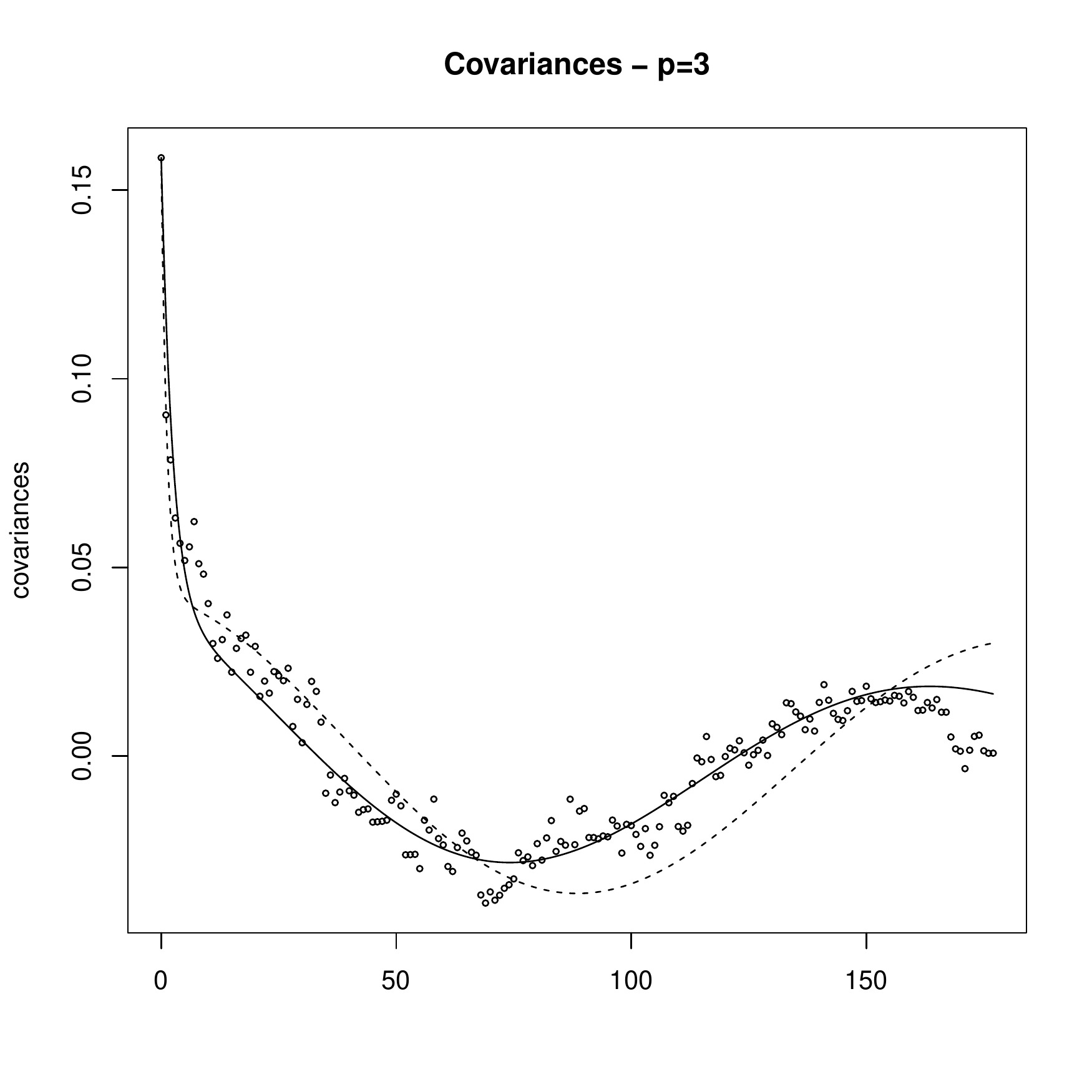}\includegraphics[width=4.5cm]{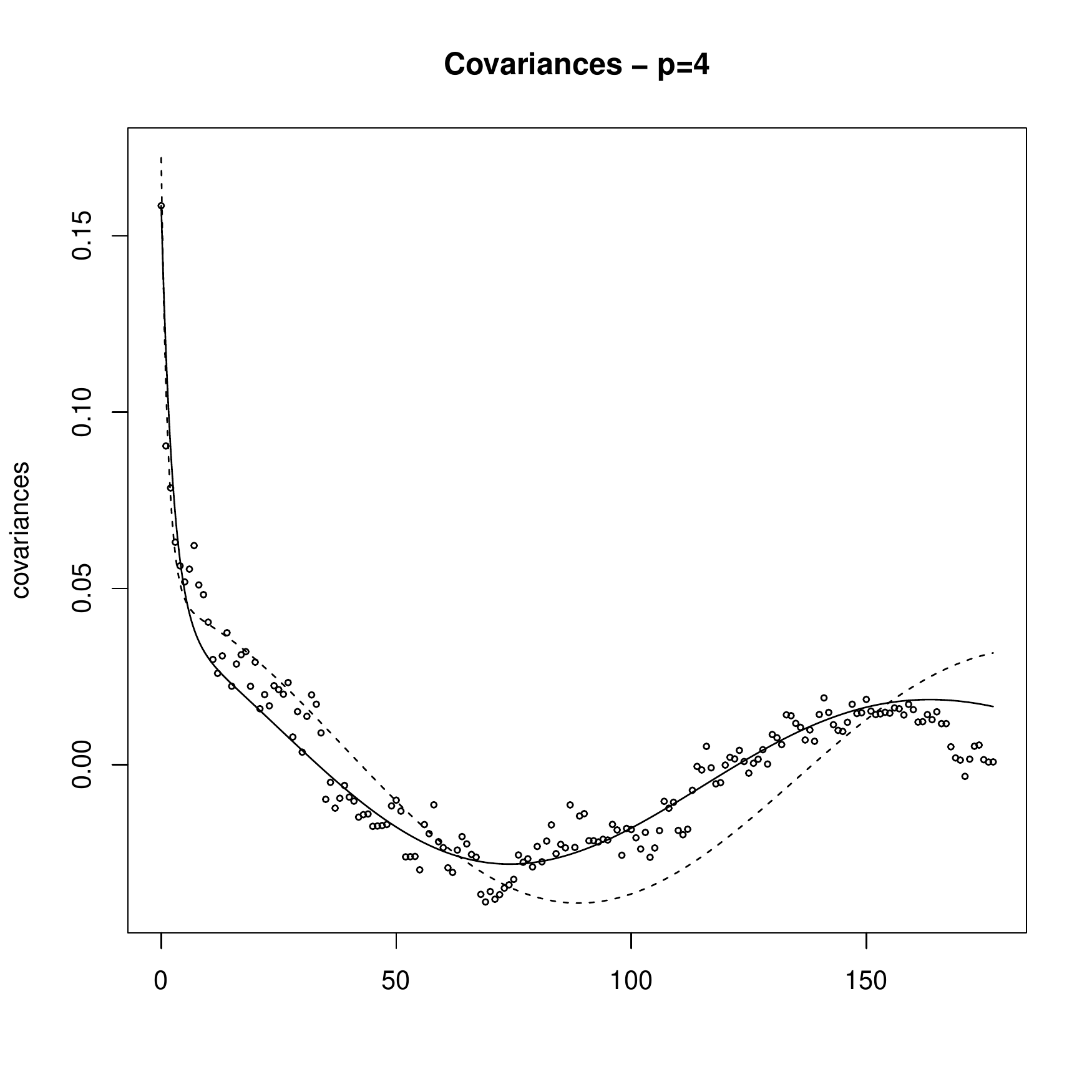}}
\caption{Empirical covariances ($\circ$) and covariances of the MC (---) and ML (- - -) fitted OU($p$) models, for $p=2,3,4$ corresponding to Series A.}\label{oupA}
\end{figure}

The ARMA(1,1) and the AR(7) fit fairly well the autocovariances for small lags, but fail to capture the structure of autocorrelations for large lags present in the series
However, the approximations obtained with the OU(3) process reflects both the short and long dependences, as shown in Figure \ref{oupA}.

\begin{figure}[!b]
\centerline{\includegraphics[width=8cm]{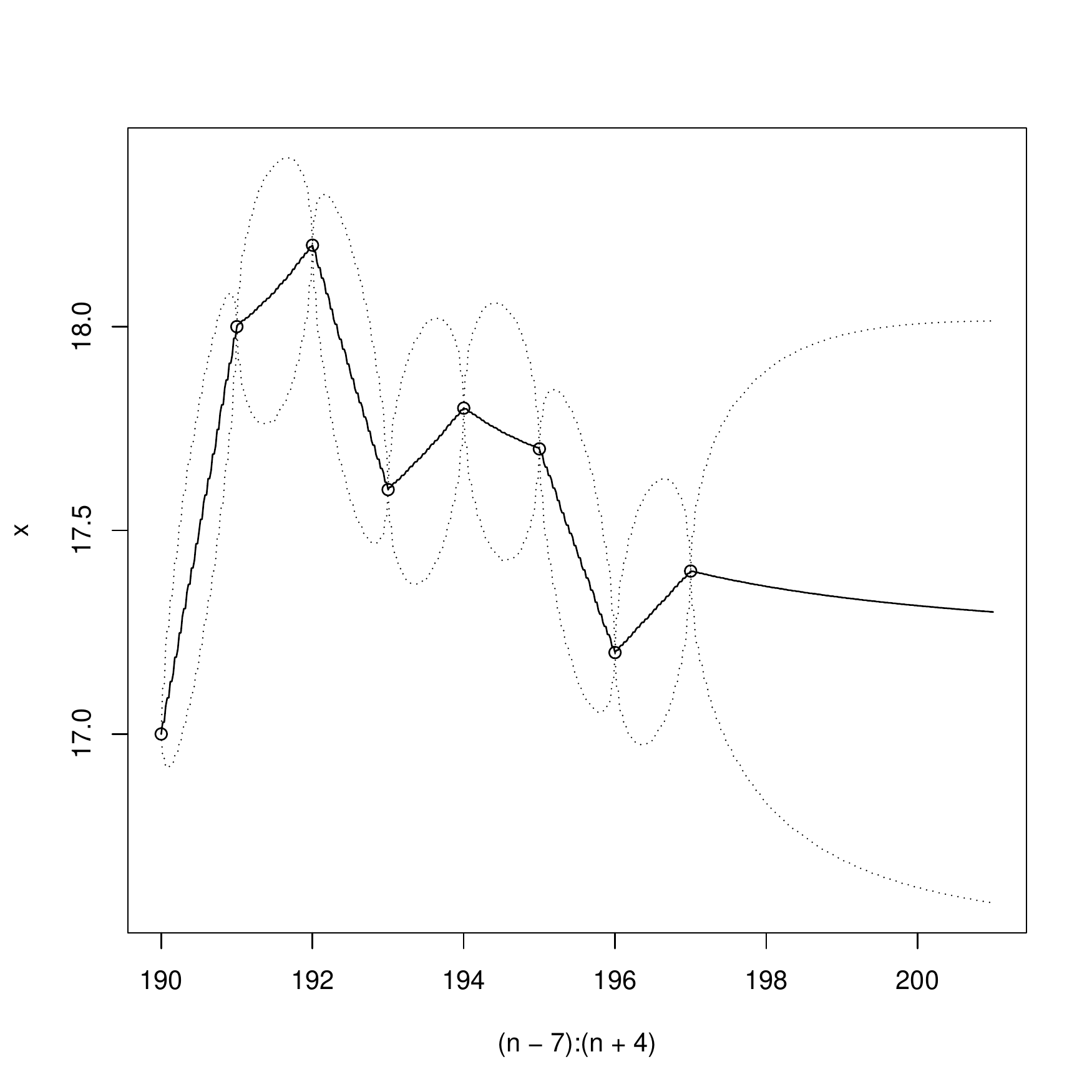}}
\caption{Confidence bands for interpolated and extrapolated values of Series A for continuous domain.} \label{SeriesA/predict}
\end{figure}

Finally we show in Figure \ref{SeriesA/predict} the predicted values of the continuous parameter process $x(t)$ for $t$ between $n-7$ and $n+4$ (190-201), obtained as the best linear predictions based on the last 90 observed values, and on the correlations given by the fitted OU(3) model. The upper and lower lines are 2$\sigma$-confidence limits for each value of the process.

\subsection{Box, Jenkins and Reinsel Series C}

The Series C is a record of $n=226$ chemical process temperature readings, taken every minute, introduced with that name in \cite{B&J}, p. 544.

\begin{figure}[!t]
\centerline{\includegraphics[width=4cm]{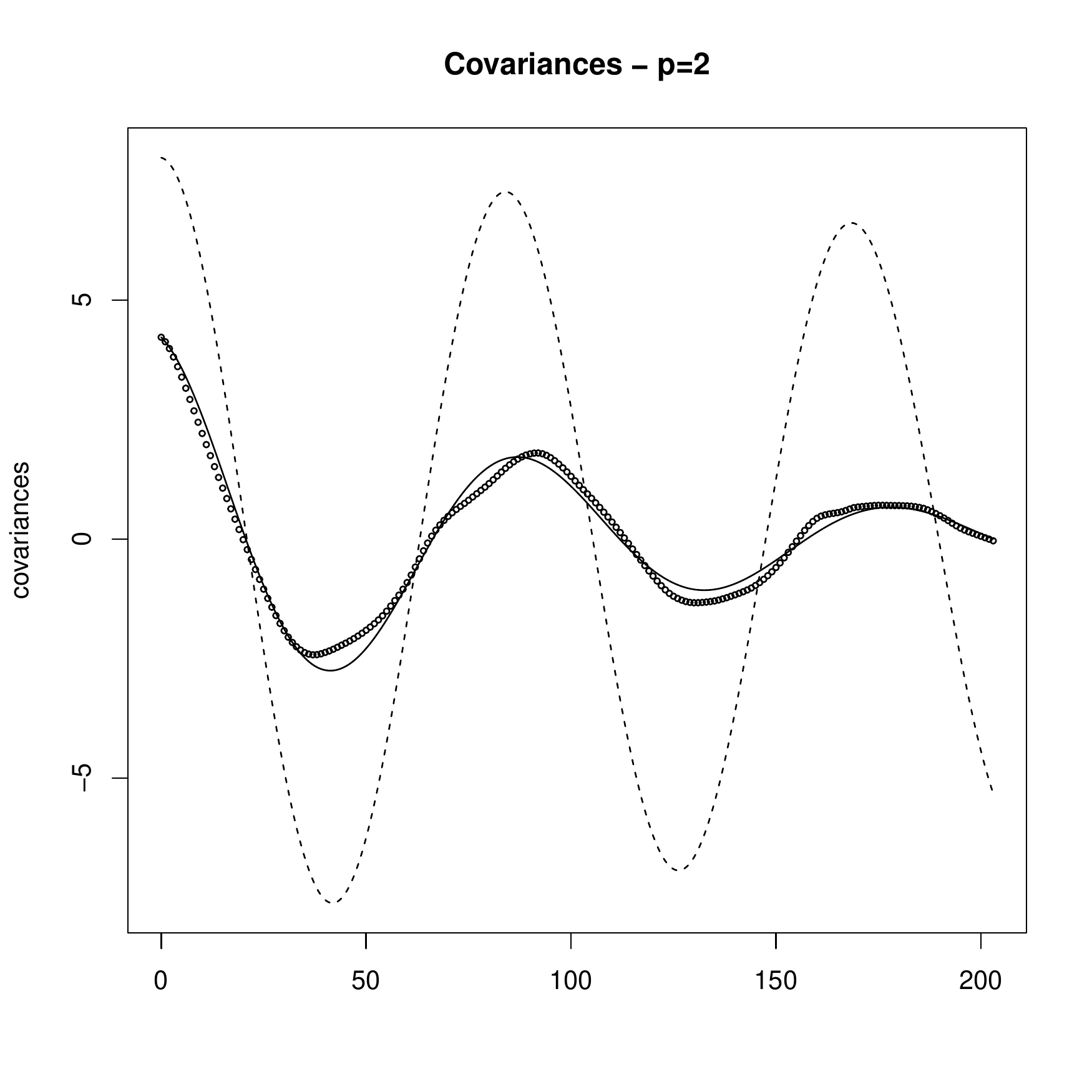}\includegraphics[width=4cm]{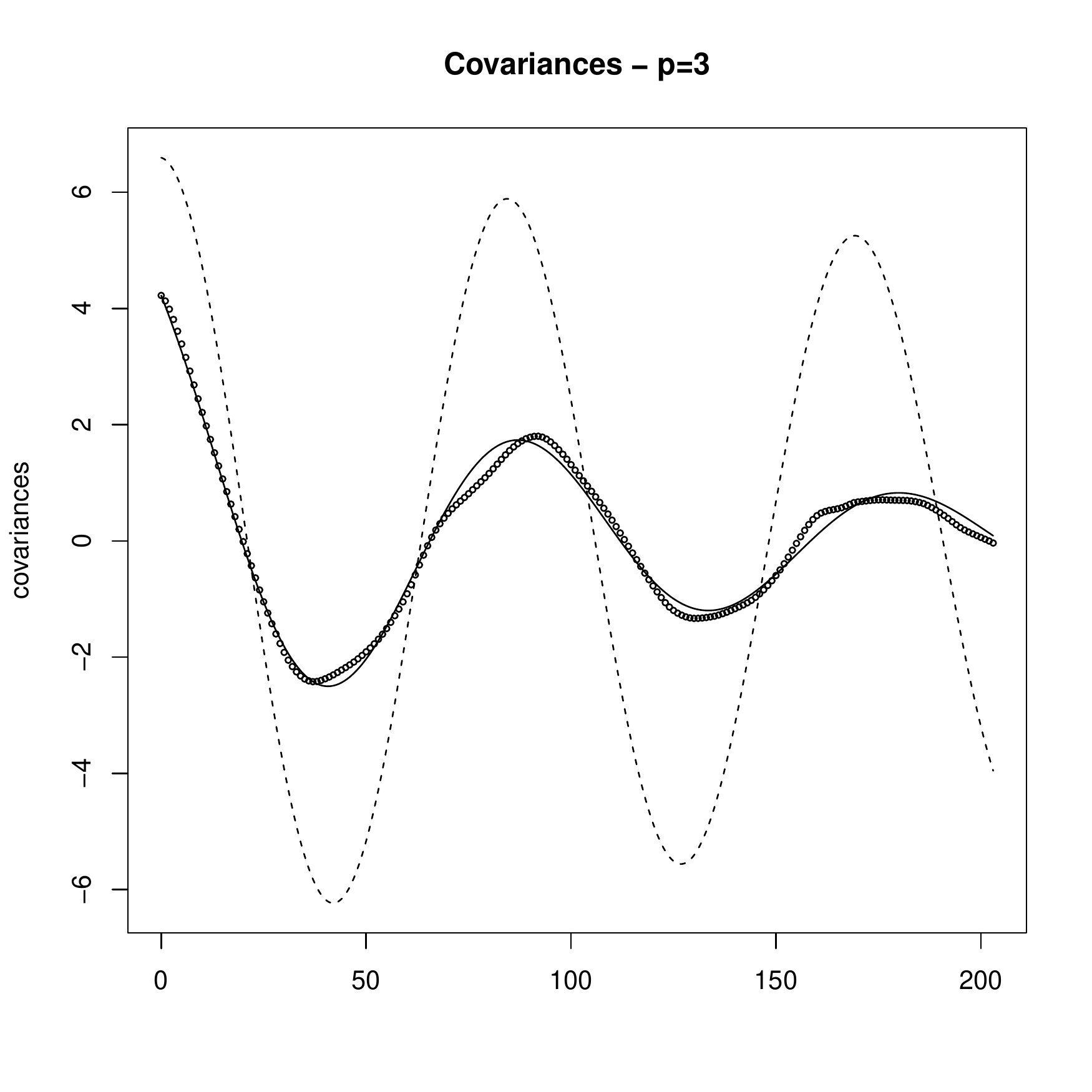}\includegraphics[width=4cm]{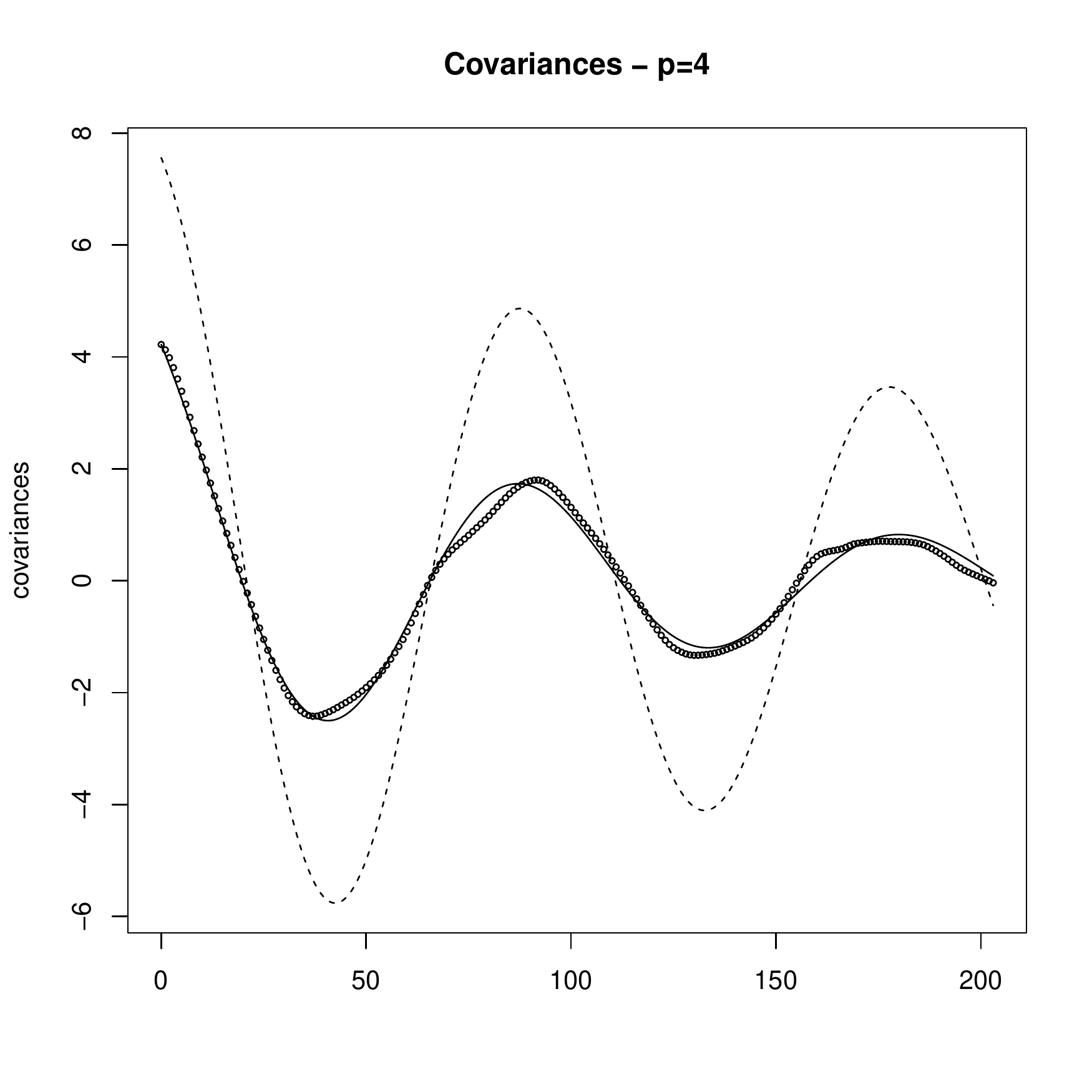}}
\caption{Empirical covariances ($\circ$) and covariances of the MC (---) and ML (- - -) fitted OU($p$) models for $p=2, 3,4$ corresponding to Series C.}\label{SeriesC/todos}
\end{figure}

As in the previous example, the fitted ARMA($p,q$) and ARIMA($p,1,q$) models for moderate values of $p$ and $q$ fail to capture the 
autocorrelations that might be present in the series.
 Figure \ref{SeriesC/todos}  shows the empirical covariances of the series and the covariances of the MC (---) and ML (- - -) fitted OU($p$) models for $p=2$, $p=3$ and $p=4$. It is not surprising that the MC estimated covariances fit better than the ML ones the empirical covariances, since they have been obtained by optimising that fit. The poor performance of the ML estimation is presumably due to the fact that the series does not obey an OU model. 

\begin{figure}[!t]
\centerline{\includegraphics[width=4cm]{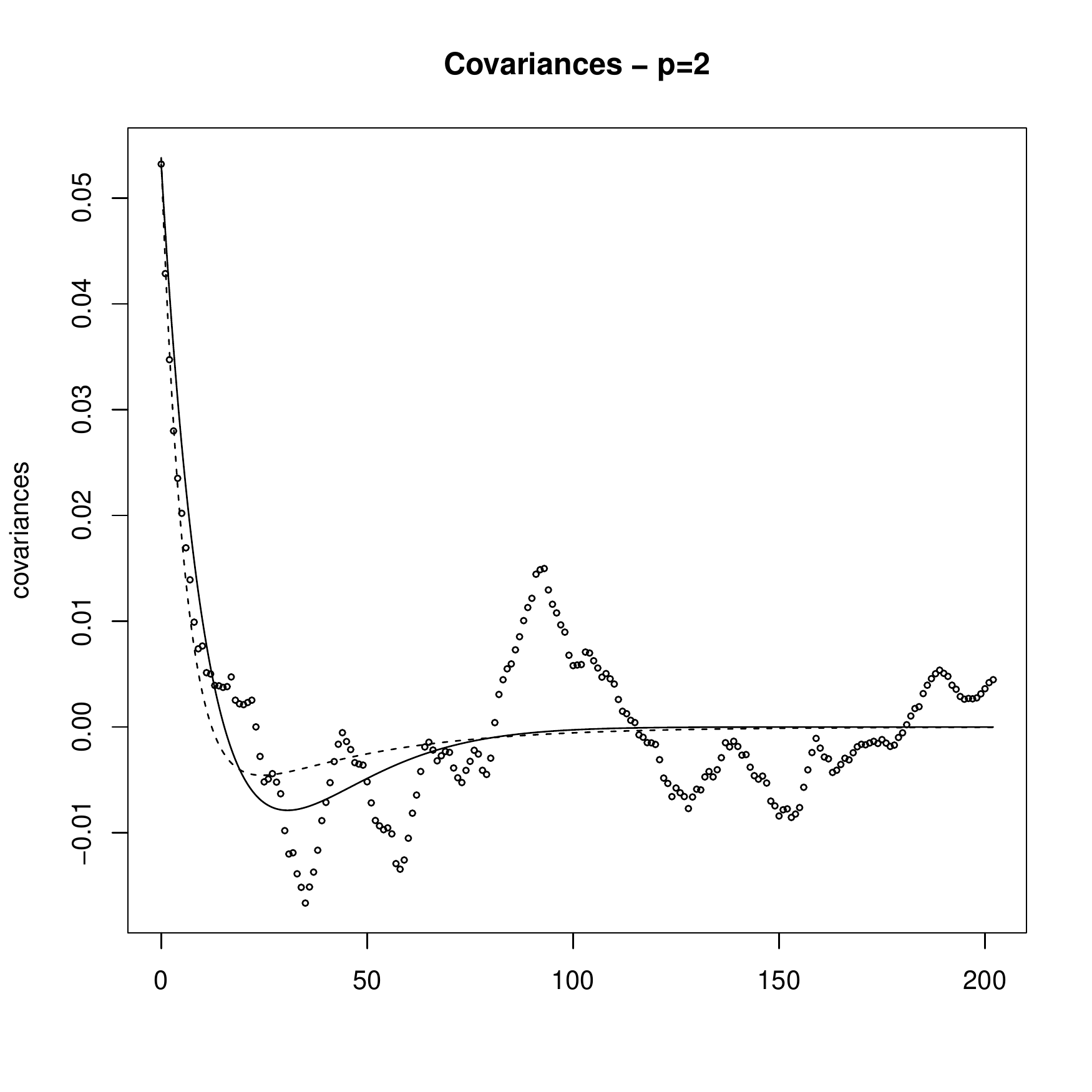}\includegraphics[width=4cm]{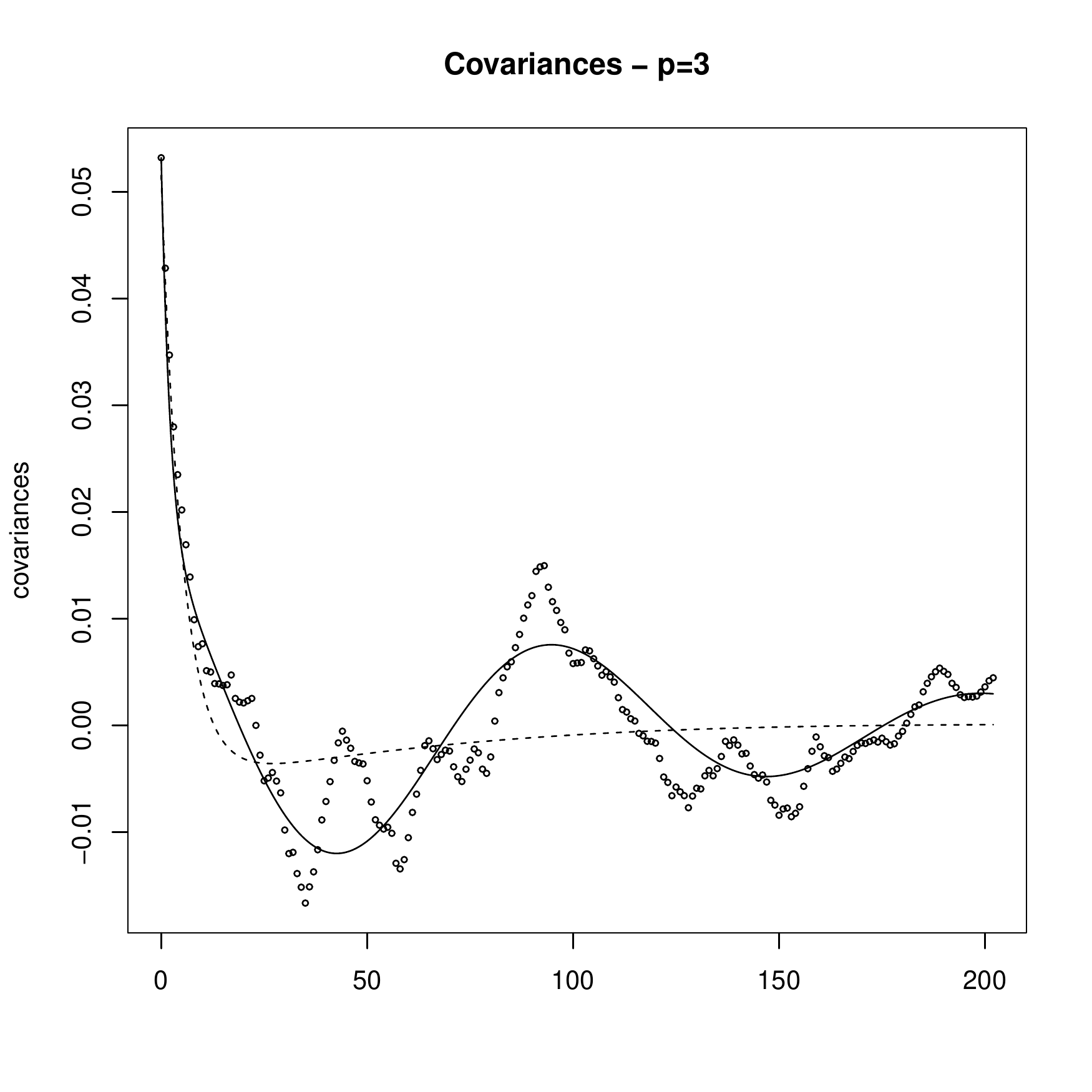}\includegraphics[width=4cm]{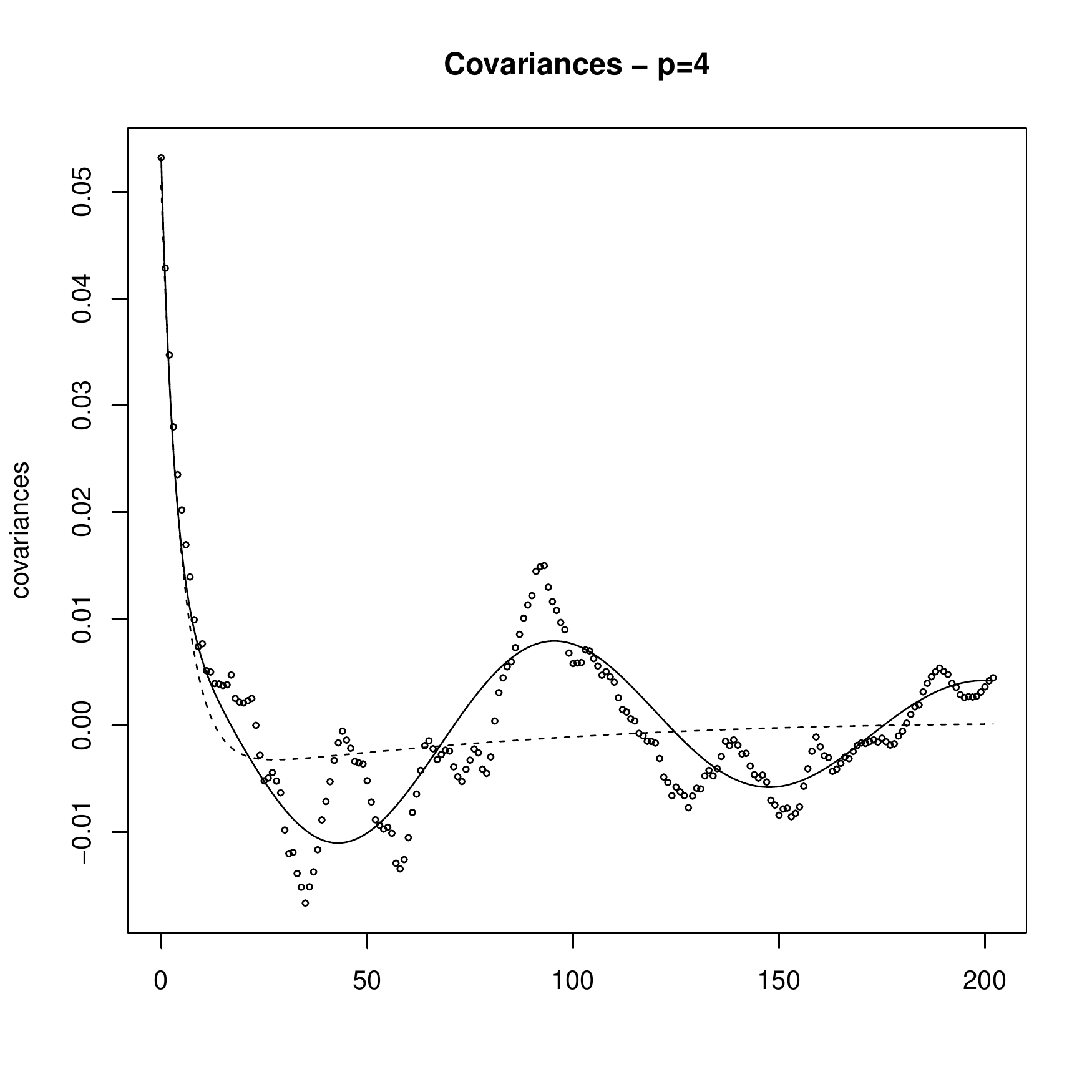}}
\caption{Empirical covariances ($\circ$) and covariances of the MC (---) and ML (- - -) fitted OU($p$) models for $p=2,3,4$ corresponding to the first differences of Series C.}\label{SeriesdC/todos}
\end{figure}

The corresponding graphs for the first differences of Series C are included in Figure \ref{SeriesdC/todos}.

\subsection{Oxigen saturation in blood}

The oxygen saturation in blood of a newborn child has been monitored during seventeen hours, and measures taken every two seconds.
We assume that a series $x_0,x_1,\dots,x_{304}$  of measures taken at intervals of 200 seconds is observed, and fit OU processes of orders $p=2,3,4$ to that series.

\begin{figure}[!t]
\centerline{\includegraphics[width=4cm]{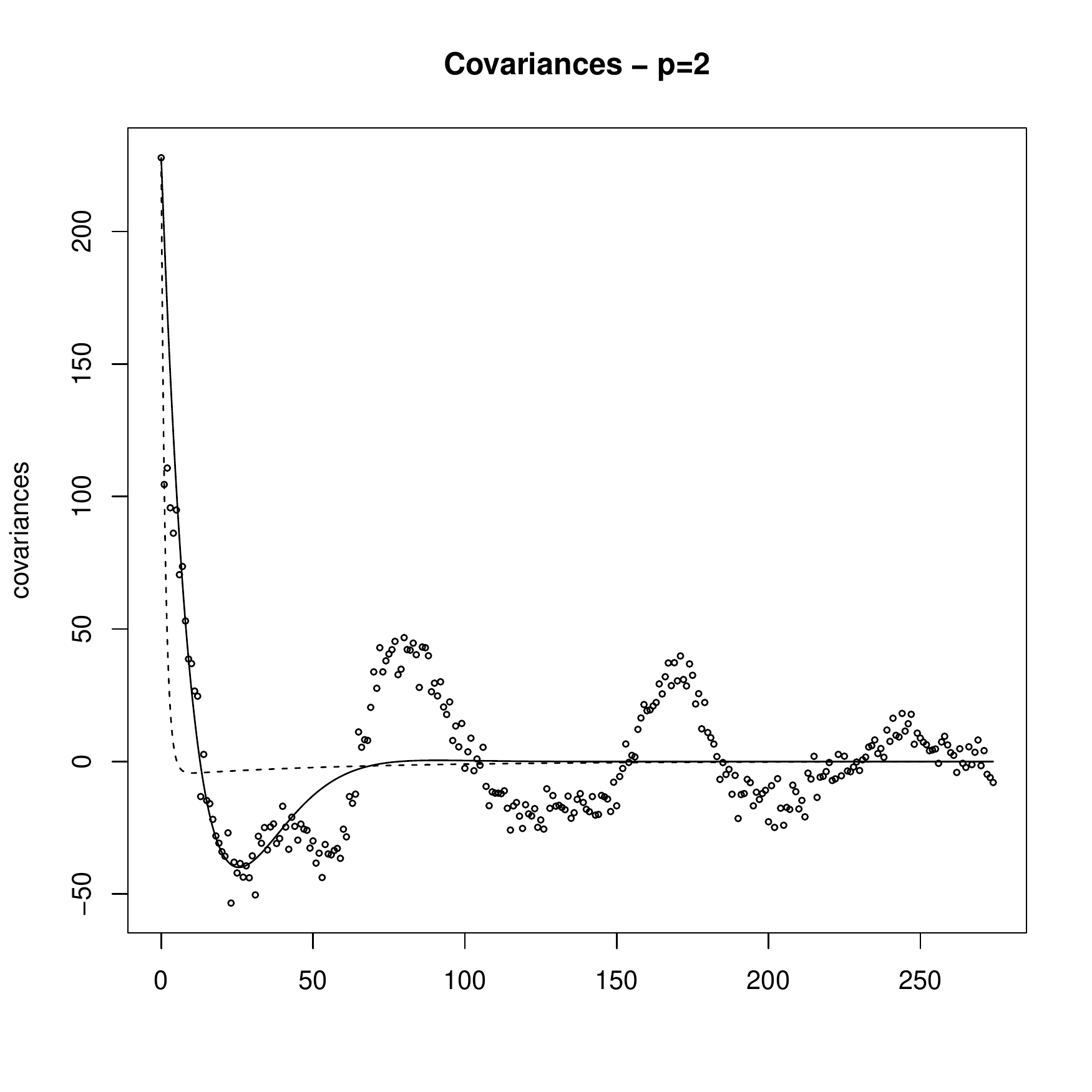}\includegraphics[width=4cm]{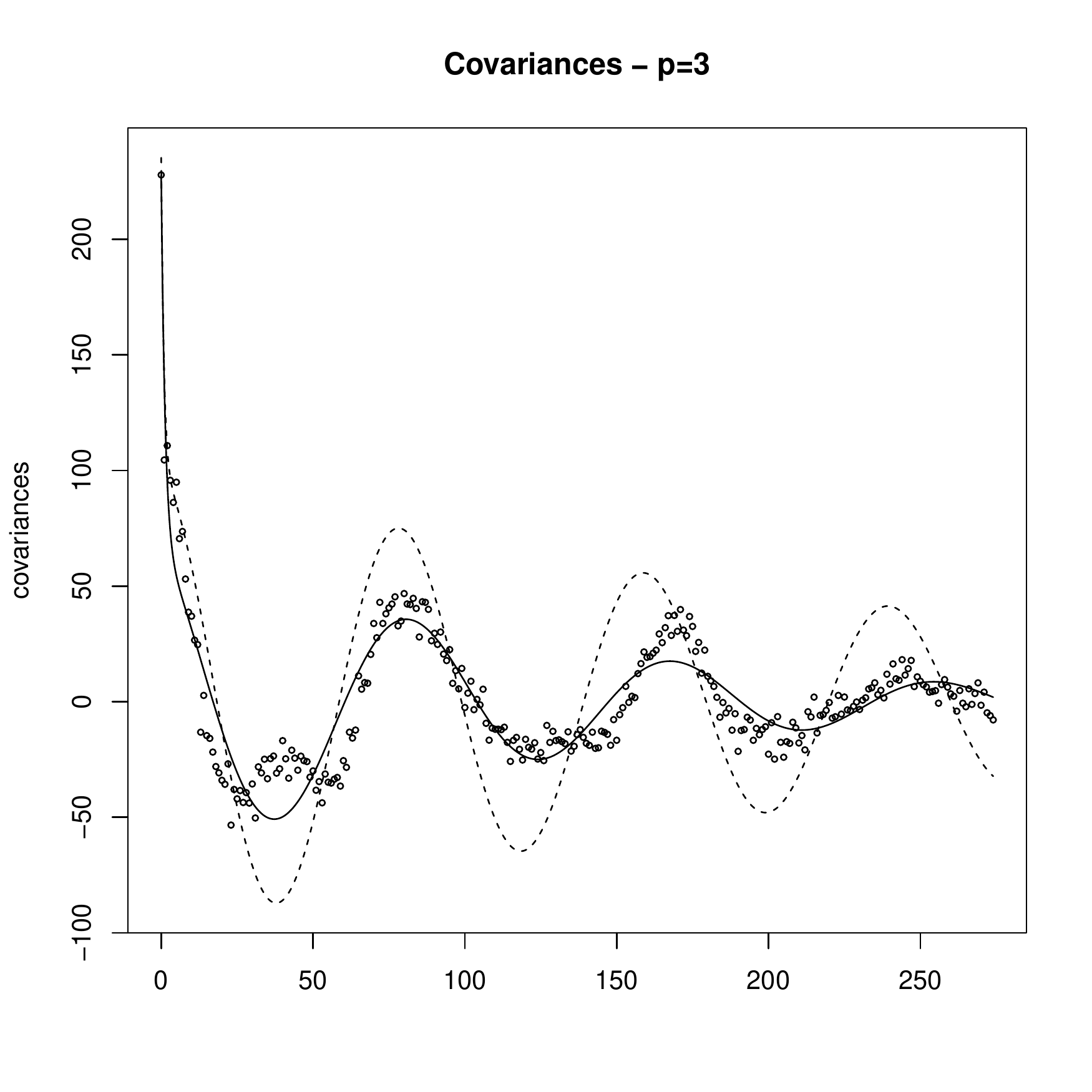}\includegraphics[width=4cm]{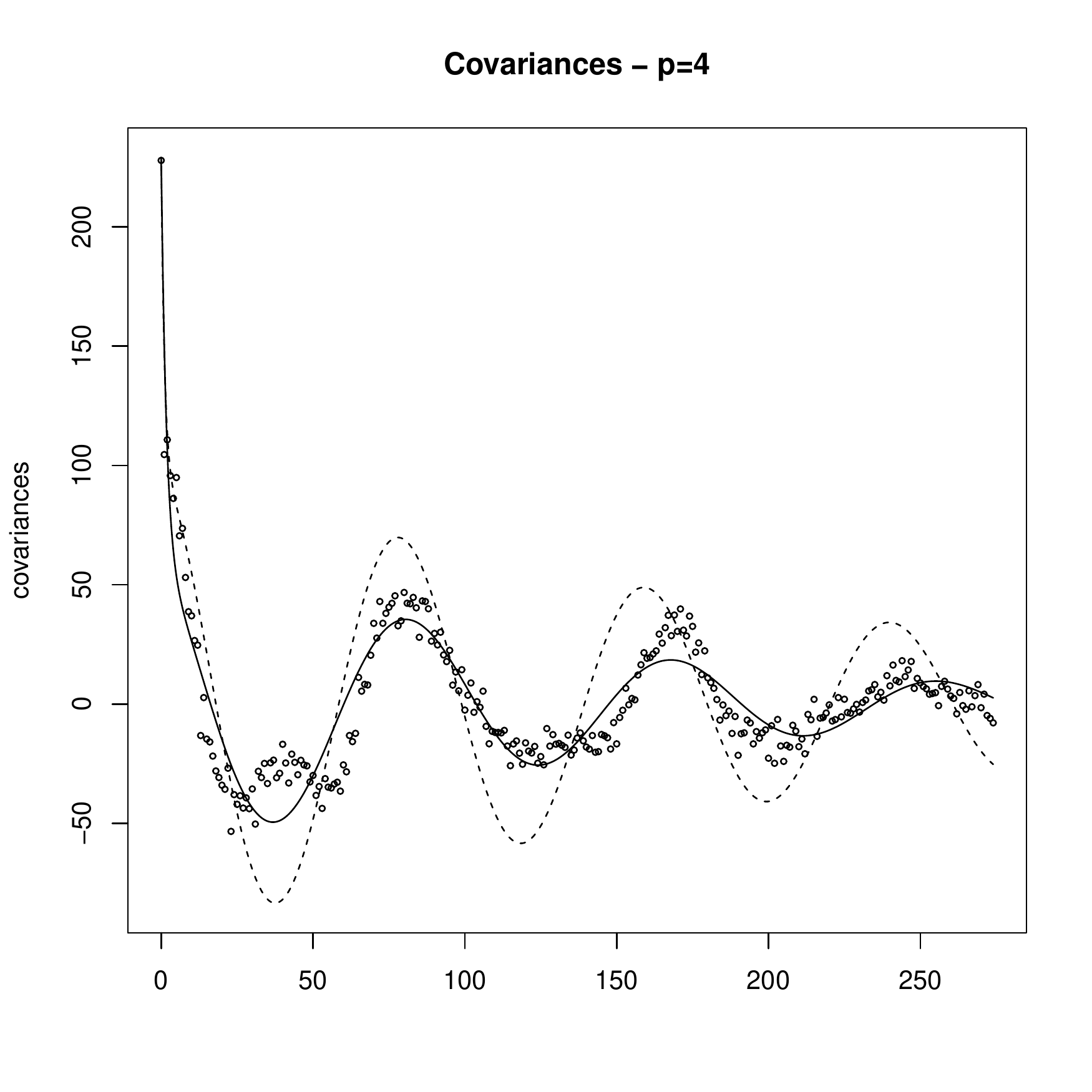}}
\caption{Empirical covariances ($\circ$) and covariances of the MC (---) and ML (- - -) fitted OU($p$) models for $p=2, 3, 4$ corresponding to the series of  O$_2$ saturation in blood.}
\label{Albanoplots}
\end{figure}

\begin{figure}[!ht]
\centerline{\includegraphics[width=8cm]{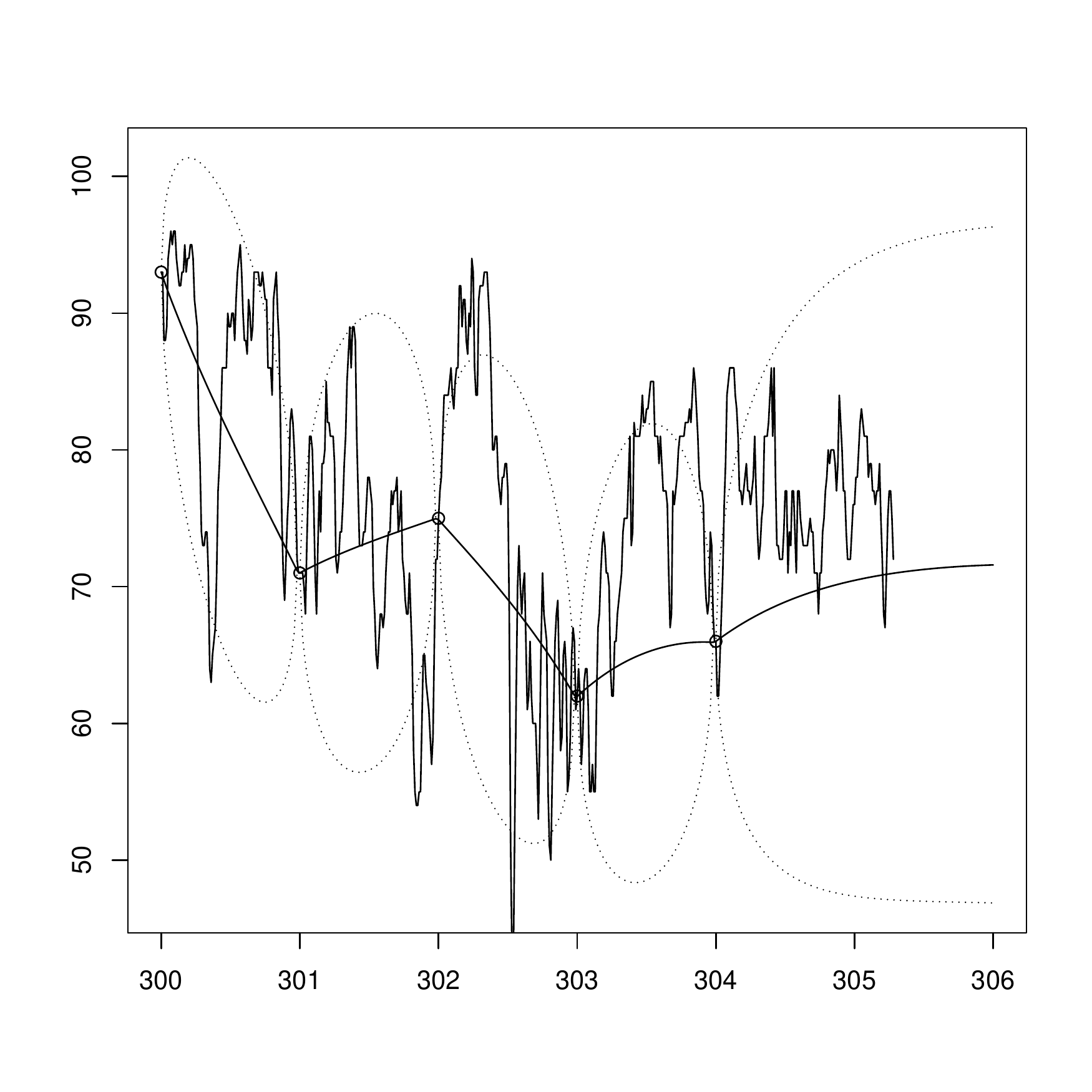}}
\caption{Partial graph showing the five last values of the series of O$_2$ saturation in blood at integer multiples of the 200 seconds unit of time ($\circ$), interpolated and extrapolated predictions (---), 2$\sigma$ confidence bands (- - -), and actual values of the series.} \label{inex}
\end{figure}

Again the empirical covariances of the series and the covariances of the fitted OU($p$) models for $p=2$, $p=3$ and $p=4$ are plotted (see Figure \ref{Albanoplots}) and the estimated interpolation and extrapolation 
are shown in Figure \ref{inex}. In the present case, the actual values of the series for integer multiples of 1/100 of the unit measure of 200 seconds are known, and plotted in the same figure.

\section{Conclusions and comments}\label{conclus}

We have proposed a family of continuous time stationary processes 
based on iterations of a generalisation of the linear operator that maps
a Wiener process onto an Ornstein-Uhlenbeck process.
The OU$(p)$ family depends on $p+1$ parameters that can be easily estimated
by either a maximum likelihood or matching correlations procedures.  Matching correlation estimators provide a fair estimation of the covariances of the data, even if the model is not well specified. 

The families  of OU$(p)$ models can be used as an alternative to ARMA or AR models for the study of 
stationary time series. For $p=1$, OU$(1)$ observed at equally spaced time instants coincides with AR$(1)$ but for larger values of $p$ the  covariances that can be described with OU$(p)$ models are not in general the same as those given by the AR$(p)$ models.  In fact, the autocorrelation structure that might be present in the data for large lags can be   modeled with OU$(p)$ with small values of $p$, a fitting that  the ARMA and AR models fail to accomplish.

\end{document}